\documentclass{amsart}

\usepackage{amsfonts,amssymb,amscd,amsmath,enumerate,verbatim,calc,xypic}
\usepackage[all]{xy} \SelectTips{eu}{}
 \usepackage[colorlinks]{hyperref}
\def\ds{\displaystyle}

\def\ov{\overline}

\def\wt{\widetilde}
\def\wh{\widehat}

\def\CA{{\mathcal A}}
\def\CB{{\mathcal B}}
\def\CC{{\mathcal C}}
\def\cD{{\mathcal D}}
\def\CE{{\mathcal E}}
\def\CF{{\mathcal F}}
\def\CG{{\mathcal G}}
\def\CI{{\mathcal I}}

\def\CL{{\mathcal L}}
\def\CK{{\mathcal K}}
\def\CM{{\mathcal M}}
\def\CN{{\mathcal N}}

\def\CS{{\mathcal S}}
\def\CT{{\mathcal T}}

\newcommand{\BN}{{\mathbb N}}
\newcommand{\BQ}{{\mathbb Q}}
\newcommand{\BZ}{{\mathbb Z}}

\newcommand{\mat}[2]{{\mathsf{M}_{#1}(#2)}}

\newcommand{\fa}{{\mathfrak a}}
\newcommand{\fm}{{\mathfrak m}}
\newcommand{\fn}{{\mathfrak n}}

\newcommand{\bsf}{{\boldsymbol f}}

\newcommand{\im}{\operatorname{Im}}
\newcommand{\id}[1]{\operatorname{id}^{#1}}

\newcommand{\sS}{\scriptstyle}

\newcommand{\ges}{\operatorname{\sS\geqslant}}
\newcommand{\les}{\operatorname{\sS\leqslant}}
\newcommand{\sges}{\operatorname{\sS>}}
\newcommand{\sles}{\operatorname{\sS<}}

\newcommand{\col}{\colon}
\newcommand{\dd}{\partial}

\newcommand{\HH}[2]{\operatorname{H}_{#1}(#2)}
\newcommand{\CH}[2]{\operatorname{H}^{#1}(#2)}
\newcommand{\cls}[1]{{\operatorname{cl}(#1)}}
\newcommand{\ZZ}[2]{\operatorname{Z}_{#1}{(#2)}}

\newcommand{\rank}{\operatorname{rank}}
\newcommand{\Ker}{\operatorname{Ker}}
\newcommand{\Coker}{\operatorname{Coker}}
\newcommand{\shift}{{\mathsf\Sigma}}
\newcommand{\susp}{{\sigma}}

\newcommand{\lra}{\longrightarrow}
\newcommand{\lla}{\longleftarrow}

\newcommand{\pd}[2]{\operatorname{pd}_{#1}{#2}}
\newcommand{\idim}[2]{\operatorname{id}_{#1}{#2}}

\newcommand{\gami}[2]{{\Gamma}^{#1}#2}

\newcommand{\gam}[1]{{\Gamma}#1}
\newcommand{\Hom}[3]{\operatorname{Hom}_{#1}({#2},{#3})}
\newcommand{\Homgr}[3]{\operatorname{Homgr}_{#1}({#2},{#3})}
\newcommand{\BigHom}[3]{\operatorname{Hom}_{#1}\!\!\bigg({\!#2\,},{#3\!\!}\bigg)}

\newcommand{\Homb}[3]{\overline{\operatorname{Hom}}_{#1}({#2},{#3})}
\newcommand{\Hombgr}[3]{\overline{\operatorname{Homgr}}_{#1}({#2},{#3})}
\newcommand{\Homv}[3]{\wh{\operatorname{Hom}}_{#1}({#2},{#3})}
\newcommand{\Homvgr}[3]{\wh{\operatorname{Homgr}}_{#1}({#2},{#3})}

\newcommand{\ExtA}{{\operatorname{Ext}}}
\newcommand{\ExtB}{\overline{\operatorname{Ext}}}
\newcommand{\ExtV}{\wh{\operatorname{Ext}}}
\newcommand{\Extb}[4]{\overline{\operatorname{Ext}}{\vphantom E}^{#1}_{#2}(#3,#4)}
\newcommand{\Extv}[4]{\wh{\operatorname{Ext}}{\vphantom E}^{#1}_{#2}(#3,#4)}

\newcommand{\BigExt}[4]{\operatorname{Ext}^{#1}_{#2}\!\!\bigg({\!#3\,},{#4\!\!}\bigg)}
\newcommand{\Ext}[4]{\operatorname{Ext}^{#1}_{#2}(#3,#4)}
\newcommand{\Tor}[4]{\operatorname{Tor}_{#1}^{#2}(#3,#4)}

\newcommand{\var}{{\hskip1pt\vert\hskip1pt}}
\newcommand{\depth}[2]{\operatorname{depth}_{#1}{#2}}
\newcommand{\edim}{\operatorname{edim}}
\newcommand{\mult}{\operatorname{mult}}
\newcommand{\codim}{\operatorname{codim}}
\newcommand{\codepth}{\operatorname{codepth}}

\newcommand{\xra}{\xrightarrow}
\newcommand{\xla}{\xleftarrow}

\newcommand{\ann}{\operatorname{Ann}}

\theoremstyle{plain}
\newtheorem{theorem}{Theorem}[section]
\newtheorem{corollary}[theorem]{Corollary}
\newtheorem{lemma}[theorem]{Lemma}
\newtheorem{sublemma}{Lemma}[subsection]

\newtheorem{subproposition}[sublemma]{Proposition}
\newtheorem{subtheorem}[sublemma]{Theorem}
\newtheorem{proposition}[theorem]{Proposition}

\theoremstyle{definition}

\newtheorem{chunk}[theorem]{}
\newtheorem{subchunk}[sublemma]{}

\newtheorem{remark}[theorem]{Remark}

\newtheorem{example}[theorem]{Example}

\numberwithin{equation}{theorem}

\theoremstyle{remark}
\newtheorem*{noremark}{Remark}
\newtheorem{question}[theorem]{Question}

\newtheorem{Remark}[equation]{Remark}

\newcommand{\numberseries}{\mdseries}   

\newlength{\thmtopspace}                
\newlength{\thmbotspace}                
\newlength{\thmheadspace}               
\newlength{\thmindent}                  

\newtheoremstyle{bfupright head,italic body}
                {\thmtopspace}{\thmbotspace}
                {\itshape}{\thmindent}{\bfseries}{}{\thmheadspace}
                {{\numberseries \thmnumber{\bf(#2)}}\thmnote{#3}}

\theoremstyle{bfupright head,italic body}
\newtheorem{ires}{}             \newtheorem*{ires*}{}

\begin{document}

\date{\today}
\title[Stable cohomology]{Stable cohomology over local rings}

\author[L.~L.~Avramov]{Luchezar~L.~Avramov}
\address{Department of Mathematics, University of Nebraska, Lincoln, Nebraska~68588}
\email{avramov@math.unl.edu}

\author[O.~Veliche]{Oana~Veliche}
\address{Department of Mathematics, University of Utah, Salt Lake City, Utah~84112}
\email{oveliche@math.utah.edu}

\thanks {L.L.A. partly supported by NSF grants DMS 0201904.}


\begin{abstract}
For a commutative noetherian ring $R$ with residue field $k$
stable cohomology modules $\Extv{n}Rkk$ have been defined for each
$n\in\BZ$, but their meaning has remained elusive.  It is proved that
the $k$-rank of any $\Extv{n}Rkk$ characterizes important properties
of $R$, such as being regular, complete intersection, or Gorenstein.
These numerical characterizations are based on results concerning the
structure of $\BZ$-graded $k$-algebra carried by stable cohomology.
It is shown that in many cases it is determined by absolute cohomology
through a canonical homomorphism of algebras $\Ext{}Rkk\to\Extv{}Rkk$.
Some techniques developed in the paper are applicable to the study of
stable cohomology functors over general associative rings.
 \end{abstract}

 \maketitle

\tableofcontents


\section*{Introduction}

A stable cohomology theory over an associative ring $R$ associates to
every pair $(M,N)$ of $R$-modules groups $\Extv nRMN$, which are defined
for each $n\in\BZ$ and vanish for all $n$ when $M$ or $N$ has finite
projective dimension.  Different constructions, have been proposed by
Benson and Carlson, Mislin, and Vogel, and all yield canonically isomorphic
theories.  However, there have been few applications outside of group
theory and Galois theory, for which the prototype---Tate cohomology for
finite groups---was created in the 1950s.

In the first four sections we develop general techniques for computing 
stable cohomology.   We approach it through a canonical transformation 
$\iota\col \ExtA_R\to\ExtV_R$ of absolute cohomology into stable
cohomology, which we study by systematically using the compatible 
multiplicative structures carried by the two theories.  A new feature are 
extensive applications of
a third cohomological functor, the bounded cohomology $\ExtB_R$, which
appears in a long exact sequence measuring the kernel and the cokernel
of $\iota$.  By extending a construction of Eisenbud we show how to track 
changes in stable cohomology under factorizations of non-zero-divisors.

The core of the paper is its second part, devoted to stable cohomology
over commutative noetherian local rings.  One goal is to investigate if
and how this theory reflects or detects properties of a ring or a module.
A second goal is to study the structure of the local cohomology functors
themselves.  Historical precedent in commutative algebra points to the
residue field $k$ of a local ring $R$ as the ultimate test case, so the
focus is kept on it for much of the second part of the paper.

When applying the general machinery to a local ring $R$ with residue
field $k$ we heavily use the fact that the absolute cohomology algebra
$\CE=\Ext{}Rkk$ is the universal enveloping algebra of a graded Lie
algebra.  The existence of such a structure underlies a well documented
successful interaction between local algebra and rational homotopy
theory. F\'elix, Halperin, and Thomas have transplanted from algebra
and systematically used in topology a notion of depth of cohomology
modules. We take the concept back into algebra and use it in a different
manner. Background material is developed in Section \ref{Depth over
graded algebras} and Section \ref{Depth of cohomology modules}.

In Section \ref{Finiteness of stable cohomology} we give necessary
and sufficient conditions for a local ring $R$ to be regular
(respectively, complete intersection, Gorenstein) in terms of the
vanishing (respectively, size, finiteness) of $\rank_k\Extv nRkk$ for
a \emph{single} value $n\in\BZ$. The last result is surprising: unlike
regularity or complete intersection, Gorensteinness is not recognized
even by the entire sequence $(\rank_k\Ext{n}Rkk)_{n\ges0}$.

In Section \ref{Structure of stable cohomology algebras} we start
the study of the graded $k$-algebra $\CS=\Extv nRkk$.  A result
of Martsinkovsky, for which we give a short proof, shows that when
$R$ is singular the map of graded algebras $\iota\col\CE\to\CS$ is
injective. We determine $\Coker(\iota)$ as a left $\CE$-module and
prove that $\depth{}\CE\ge2$ implies $\CS=\iota(\CE)\oplus\CT$ where
$\CT$ is the $\CE$-torsion submodule of $\CS$, and is the \emph{unique}
direct complement of $\iota(\CE)$ as a left $\CE$-module.  In Section
\ref{Stable cohomology algebras of complete intersection rings} this
information is used to produce a nearly complete, explicit computation
of the algebra $\CS$ for complete intersection rings.

It is natural to ask whether the results for complete intersections
extend, in some form, to all singular Gorenstein rings $R$ with $\codim
R\ge2$.  In Section \ref{Stable cohomology over Gorenstein rings}
we prove that $\depth{}\CE\ge2$ implies $\CT$ is a two-sided ideal of
$\CS$ with $\CT^2=0$, is isomorphic to a shift of $\Hom k\CE k$ as left
$\CE$-module, and $\CS=\iota(\CE)\oplus\CT$.  A similar relationship
between the Tate cohomology algebra $\wh{\rm H}{}^*(G,k)$ of a finite
group $G$ and its cohomology algebra $\CH{*}{G,k}$ was discovered by
Benson and Carlson. The parallel is remarkable, as $\CH{*}{G,k}$ is
graded-commutative and finitely generated, while $\CE$ may be non-finitely
generated and almost always is very far from commutative.

One has $\depth{}\CE\ge1$ for all singular rings, so the condition
$\depth{}\CE\ge2$ is not too special.  We prove that it holds for
several classes of Gorenstein rings, including those of codimension $2$
or $3$, those of minimal multiplicity, and the localizations of Koszul
algebras. We are not aware whether splitting occurs always: Gorenstein
rings with $\depth{}\CE=1$ are hard to come by, and for the known ones
$\CE$ splits off $\CS$.

It was noted above that when $R$ is not Gorenstein $\rank_k\Extv nRkk$ is
infinite for each $n$, so over such rings a different structure of $\CS$
may be expected.  As a test case in Section \ref{Stable cohomology over
Golod rings} we turn to Golod rings, whose homological properties
are in many respects opposite to those of Gorenstein rings.  We show
that $\depth{}\CE=1$ holds for all Golod rings, and for a subclass of
such rings we work out the structure of $\CS$ in sufficient detail to
prove that $\iota$ does not split as a map of left $\CE$-module.

\section{Cohomology theories}
\label{Cohomology theories}

We start by describing notions concerning complexes and, more generally, DG
(that is, differential graded) modules and algebras.  The latter are used to describe
composition products carried by the \emph{absolute cohomology} functors. 
We then introduce a \emph{bounded cohomology} theory that has not been systematically studied before.  Finally, we present Vogel's construction 
of \emph{stable cohomology}.

\begin{subsection}{DG algebras and DG modules}
 \label{shift}
To grade a complex $C$ we use subscripts \emph{or} superscripts.  Thus,
$C$ can be written either as a sequence of maps $\dd_n^C\col C_n\to
C_{n-1}$, or as a sequence of maps $\dd^{-n}_C\col C^{-n}\to C^{-n+1}$,
with $\dd^{-n}_C=\dd^C_n$. Accordingly, an element $c\in C_n$ is assigned
a lower (or homological) degree $n$, denoted \emph{and} an upper (or
cohomological) degree $-n$; we write $\lfloor c\rfloor=n$ and $\lceil
c\rceil=-n$, respectively.  When the nature of degree does not matter
we use $|c|$ in place of either $\lfloor c\rfloor$ or $\lceil c\rceil$.

When $z\in C$ is a cycle $\cls{z}$ denotes its homology class.  

For every $s\in\BZ$ let $\shift^sC$ denote the complex with
$(\shift^sC)_n=C_{n-s}$ and $\dd^{\shift^sC}_n=
(-1)^s\dd^{C}_{n-s}$; let $\susp^s\col C\to\shift^sC$ be the
bijective map $C_n\ni c\mapsto c\in(\shift^sC)_{n+s}$.

Bimodules have actions from the left and from the right,
listed in that order.  If $A$ and $A'$ are DG algebras and $C$ is
a DG $A$-$A'$-bimodule, then the formula
 \[
a\cdot\susp^s(c)\cdot a'=(-1)^{|a|s}\susp^s(a\cdot c\cdot a')
 \]
turns $\shift^s C$ into a DG $A$-$A'$-bimodule and $\susp^s$
into a chain map of DG bimodules. Furthermore, when $B$ is a DG
$A$-$A$-bimodule the map
\begin{align*}
B\otimes_A\shift^sC\lra\shift^s(B\otimes_AC)
\end{align*}
given by $b\otimes\susp^s(c)\mapsto(-1)^{|b|s}\susp^{s}(b\otimes
c)$ is an isomorphism of DG $A$-$A'$-modules.
 \end{subsection}
 \bigskip

For the rest of the section $R$ denotes an associative ring,
$M$ and $N$ are left $R$-modules, and
and $F\to M$ and $G\to N$ are projective resolutions.  

\begin{subsection}{Absolute Ext}
 \label{absolute} 
Let $\Hom RFG$ denote the complex of abelian groups with
\[
\Hom RFG_n=\prod_{i\in\BZ}\Hom R{F_i}{G_{i+n}}=\Hom RFG^{-n}
\]
as component of homological degree $n$ (cohomological degree
$-n$), and differential
\[
\dd(\beta)=\dd^G\beta-(-1)^{|\beta|}\beta\dd^F\,.
\]
The induced map $\Hom RFG\to\Hom RFN$ is a quasi-isomorphism, so
one has
\[
\CH{}{\Hom RFG}\cong\CH{}{\Hom RFN}=\Ext{}RMN\,.
\]

\begin{subchunk}
 \label{hom}
Composition of homomorphisms turns $\Hom RFF$ and $\Hom RGG$ into
DG algebras, and $\Hom RFG$ into a DG $\Hom RGG$-$\Hom
RFF$-bimodule.  The \emph{composition products} induced in
homology can be computed from any pair of projective resolutions.
They turn $\Ext{}RMM$ and $\Ext{}RNN$ into graded algebras, and
$\Ext{}RMN$ into a graded $\Ext{}RNN$-$\Ext{}RMM$-bimodule.
 \end{subchunk}

\begin{subchunk}
 \label{tensor}
The DG algebra $\Hom RGG$ acts on the complex $G$ by evaluation of
homomorphisms.  For every complex $C$ of right $R$-modules the map
\begin{align*}
\Hom RGG\otimes_{\BZ}(C\otimes_RG)&\lra C\otimes_RG\\
\alpha\otimes(c\otimes g)&\longmapsto(-1)^{|\alpha||c|}c\otimes
\alpha(g)
\end{align*}
endows $C\otimes_RG$ with a structure of left DG module over $\Hom
RGG$. Clearly, this structure is natural with respect to morphisms
of complexes $C\to C'$.
 \end{subchunk}

\begin{subchunk}
 \label{tor}
\label{Tor} Let $L$ be a right $R$-module.  Setting $C=L$ in
\eqref{tensor} one obtains a morphism
\[
\Hom RGG\otimes_{\BZ}(L\otimes_RG)\lra L\otimes_RG\,.
\]
In homology it induces for all $l,n\in\BZ$ homomorphisms of abelian
groups
\[
\Ext nRNN\otimes_\BZ\Tor lRLN\lra\Tor{l-n}RLN
\]
that turn $\Tor{}RLN$ into a graded left module over $\Ext{}RNN$.
 \end{subchunk}
 \end{subsection}

\begin{subsection}{Bounded Ext}
\label{bounded} 
A homomorphism $\beta\in\Hom RFG$ is
\emph{bounded} if $\beta_i=0$ for all $i\gg0$.  The subset $\Homb
RFG$ of $\Hom RFG$, consisting of all bounded homomorphisms, is 
a subcomplex, with components
 \[
\Homb RFG_n=\coprod_{i\in\BZ}\Hom R{F_{i}}{G_{i+n}}=\Homb
RFG^{-n}\,.
 \]
The graded abelian group $\Extb{}RMN=\CH{}{\Homb RFG}$
with components
\[
\Extb{n}RMN=\CH{n}{\Homb RFG}\,,
\]
is called the \emph{bounded cohomology} of $M$ and $N$ over $R$. 

\begin{subchunk}
 \label{bounded products}
It is easy to see that $\Homb RFG$ is a DG subbimodule of $\Hom
RFG$ for the actions of $\Hom RFF$ and $\Hom RGG$ described in
\eqref{hom}, so $\Extb{}RMN$ becomes a graded
$\Ext{}RNN$-$\Ext{}RMM$-bimodule. 
  \end{subchunk}

The elementary observation below plays a pivotal role in the
paper.  It should be noted that the right-hand analog of this
statement fails, see Example \eqref{non-nilpotent}.

\begin{sublemma}
 \label{nilpotent}
For every $\tau\!\in\!\Extb{}RMN$ there exists an integer $j\ge0$,
such that
 \[
\Ext{\ges j}RNN\cdot\tau=0\,.
 \]
 \end{sublemma}

\begin{proof}
By hypothesis, $\tau=\cls\beta$ for some chain map $\beta\in\Hom RFG$ satisfying
$\beta(F)\subseteq G_{\sles j}$ for some $j\geq 0$.
For each $\gamma\in\Hom R GG_n$ one then has
\[
(\gamma\beta)(F)=\gamma\big(\beta(F)\big)
\subseteq\gamma\big(G_{\sles j}\big)\subseteq G_{\sles j+n}
\]
Since $G_{\sles j+n}=0$ for $n\le-j$, this implies
$\Ext{\ges j}RNN\cdot\cls\beta=0$.
 \end{proof}

Some of the DG module structures discussed so far are linked as
follows:

\begin{sublemma}
\label{omega} There is a morphism of DG $\Hom RGG$-$\Hom
RFF$-bimodules
\[
\omega\col\Hom RFR\otimes_RG\lra\Homb RFG
\]
with actions on the source given by \eqref{hom}, \eqref{tensor},
and on the target by \eqref{bounded products}.

If the $R$-module $F_i$ is finite for each $i$, then $\omega$ is
bijective.
 \end{sublemma}

\begin{proof}
It is easy to verify that \( \omega'(\phi\otimes
g)(f)=(-1)^{|g||f|}\phi(f)g \) defines a morphism
 \[
\omega'\col \Hom RFR\otimes_RG\lra\Hom RFG
 \]
of DG bimodules.  The image of $\omega'$ lies in $\Homb RFG$, so
it yields a morphism $\omega$ with the desired source and target.
For each $n\in\BZ$ the following equalities
\begin{align*}
(\Hom RFR\otimes_RG)_n&=\coprod_{(-i)+j=n}\Hom R{F_i}R\otimes_R G_{j}\\
\Homb RFG_n&=\coprod_{j-i=n}\Hom R{F_{i}}{G_j}
\end{align*}
hold by definition.  When each $R$-module $F_i$ is finite,
$\omega$ restricts to an isomorphism $\Hom R{F_i}R\otimes_R
G_{j}\to\Hom R{F_{i}}{G_j}$ for each pair $(i,j)$.
 \end{proof}
  \end{subsection}

\begin{subsection}{Stable Ext}
 \label{stable}
Using the subcomplex $\Homb RFG$ described in \eqref{bounded},
set
\[
\Homv RFG=\Hom RFG/\Homb RFG\,.
\]
Following Pierre Vogel, we define the \emph{stable cohomology} of
$M$ and $N$ over $R$ to be the graded abelian group
$\Extv{}RMN=\CH{}{\Homv RFG}$ with components
\[
\Extv{n}RMN=\CH{n}{\Homv RFG}\,.
\]
The assignment $(M,N)\mapsto\Extv{}RMN$ yields a cohomological
functor, contravariant in $M$ and covariant in $N$, from
$R$-modules to graded $\BZ$-modules.

\begin{subchunk}
 \label{stable-products}
As $\Homb RFG$ is a DG subbimodule of $\Hom RFG$ for the
left action of $\Hom RGG$ and the right action of $\Hom RFF$,
see \eqref{bounded products}, one sees that $\Extv{}RMM$ and 
$\Extv{}RNN$ are graded algebras, and
$\Extv{}RMN$ is a graded $\Extv{}RNN$-$\Extv{}RMM$-bimodule.
 \end{subchunk}

Stable cohomology over general associative rings took a long time 
to emerge, and then it appeared in several avatars.   We give a short, 
incomplete list of sources.

 \begin{subchunk}
Historically the first example of stable cohomology is Tate's 
cohomology theory $\wh{\rm H}{}^n(G,-)$ for modules over a 
finite group $G$: One has $\wh{\rm H}{}^n(G,-)=
\Extv n{\BZ G}{\BZ}-$, where $\BZ G$ is the group ring of $G$;
see \cite[Ch.\ XII]{CE}.  Tate's construction is based on complete resolutions of
$\BZ$.  Buchweitz \cite{Bu} extended the technique to define
a two-variable theory over two-sided noetherian Gorenstein rings.

The functors $\Extv nR--$ were introduced by Vogel in the
mid-1980s.  The first published account appears only in \cite{Go},
where it is called `Tate-Vogel cohomology'.  Different approaches were
independently proposed by Benson and Carlson \cite{BC} and by Mislin
\cite{Mi}; background and details can be found in Kropholler's survey
\cite[\S 4]{Kr}\footnote{Where Definition (4.2.2) contains a typo:
$\Omega^{i+n}N$ should be changed to $\Omega^{i-n}N$.}.  We have settled
on the name `stable cohomology' to emphasize the fact that $\Extv 0RMN$
is a group of homomorphisms of objects in a \emph{stabilization} of the
category of $R$-modules, see Beligiannis \cite[\S\S 3,5]{Bl} for details.
 \end{subchunk}
 
 \section{Comparisons}
 \label{comparisons}

In this section $R$ is an associative ring, $M$ and $N$ are left
$R$-modules, $F\to M$ and $G\to N$ denote projective resolutions.  
The objective is to describe important links between the
cohomology theories introduced in Section \ref{Cohomology theories}.

 \begin{chunk}
 \label{sequence}
By construction, there is an exact sequence of DG bimodules
 \begin{equation}
 \label{sequence:complexes}
0\lra\Homb RFG\lra\Hom RFG\lra\Homv RFG\lra0
 \end{equation}
that is unique up to homotopy.  It defines an exact
sequence 
 \begin{equation}
 \label{sequence:homology1}
\begin{gathered}
 \xymatrixrowsep{.7pc} \xymatrixcolsep{2pc}\xymatrix{
\Extb{}RMN\ar@{->}[r]^-{\eta}
&\Ext{}RMN\ar@{->}[r]^-{\iota}
&\Extv{}RMN\ar@{->}[r]^-{\eth}
&{}
\\
\shift\Extb{}RMN\ar@{->}[r]^-{\shift\eta} &\shift\Ext{}RMN}
 \end{gathered}
 \end{equation}
of graded $\Ext{}RNN$-$\Ext{}RMM$-bimodules.  Thus,
there is an exact sequence
 \begin{equation}
 \label{sequence:homology2}
\begin{gathered}
 \xymatrixrowsep{.7pc} \xymatrixcolsep{1.8pc}\xymatrix{
\cdots\quad\ar@{->}[r]
&\Extb{n}RMN\ar@{->}[r]^-{\eta^n}
&\Ext{n}RMN\ar@{->}[r]^-{\iota^n}
&\Extv{n}RMN\ar@{->}[r]^-{\eth^n}
&{}
\\
&\Extb{n+1}RMN\ar@{->}[r]^-{\eta^{n+1}}
&\Ext{n+1}RMN\ar@{->}[r]
&\quad\cdots
}
 \end{gathered}
 \end{equation}
of abelian groups, and the latter is natural in both module arguments.

Furthermore, $\iota\col\Ext{}RMM\to\Extv{}RMM$ is a homomorphism 
of graded algebras, and $\iota\col\Ext{}RMN\to\Extv{}RMN$ is an
equivariant homomorphism of graded $\Extv{}RNN$-$\Extv{}RMM$-bimodules.
 \end{chunk}

It is easy to determine if $\eta$ is an isomorphism; see also \cite[(4.2.4)]{Kr}, 
\cite[(4.5.1)]{Ve}.

\begin{proposition}
 \label{finite-pd}
The following conditions are equivalent.
 \begin{enumerate}[\quad\rm(i)]
  \item
$M$ has finite projective dimension.
  \item
$\Extv nR{M}{-}=0$ for every $n\in\BZ$.
  \item
$\Extv nR{-}{M}=0$ for every $n\in\BZ$.
  \item
$\Extv 0R{M}{M}=0$.
 \item
$\eta^n\col\Extb{n}RM-\to\Ext{n}RM-$ is an isomorphism for every $n\in\BZ$.
 \item
$\eta^n\col\Extb{n}R-M\to\Ext{n}R-M$ is an isomorphism for every $n\in\BZ$..
 \end{enumerate}
 \end{proposition}

 \begin{proof}
Choosing a finite resolution $F\to M$ one gets $\Homb RFG=\Hom RFG$
for every resolution $G$.   Thus, (i) implies (v) and (vi).  The 
exact sequence \eqref{sequence:homology2} shows that (v) implies (ii),
and (vi) implies (iii).  It is clear that (ii) or (iii) implies (iv).  If (iv) holds, then 
for some $\gamma\in{\Hom RFF}^{1}$ and some $p\ge0$ the morphism 
$\beta={\id F}-\dd\gamma+\gamma\dd\col F\to F$ satisfies $\beta_i=0$ for all 
$i\ge p$.  Thus, one gets
 \[
0=\CH{p}{\Hom R{\beta}-}=\CH p{\Hom R{\id F}-}=\id{\CH p{\Hom R{\beta}-}}
 \]
hence $\Ext pRM-=\CH p{\Hom RF-}=0$; that is, $\pd RM<p$, so (i) holds.
 \end{proof}

\setcounter{theorem}{1}

Next we give a criterion for $\iota$ to be an isomorphism in 
high degrees.
  
\begin{theorem}
 \label{perpendicular}
For an integer $m$ the following conditions are equivalent.
\begin{enumerate}[\rm\quad(i)]
 \item
$\iota^n\col\Ext nR{M}{-}\to\Extv nRM-$ is an isomorphism for
all $n>m$ and $\iota^{m}$ is an epimorphism.
 \item
$\Ext nR{M}{P}=0$ for all $n>m$ and every projective
$R$-module $P$.
 \item
$\Extb nR{M}{-}=0$ for all $n>m$.
 \end{enumerate}
When $M$ has a resolution by finite projective modules they are
also equivalent to
\begin{enumerate}[\rm\quad(i)]
 \item[\rm(ii$'$)]
$\Ext nRMR=0$ for all $n>m$.
\end{enumerate}
\end{theorem}

\begin{proof}
The exact sequence \eqref{sequence:homology2} shows that (i) and (iii) are
equivalent. If (iii) holds, then so does (ii), because 
$\Ext{n}R{M}{-}\cong\Extb{n}R{M}{-}$ by Proposition \eqref{finite-pd}. It
is clear that (ii) implies (ii$'$). The converse holds because the
hypothesis on $M$ in (ii$'$) implies that the functor $\Ext nRM-$
commutes with all direct sums.

To prove that (ii) implies (iii) fix an integer $n>m$ and choose a chain map
$\alpha\in\Homb RFG^n$.  Thus, for some fixed $s\ge
n$ and all $j\ge s$ one has $\alpha_j=0$, while
 \begin{equation}
\tag*{(\ref{perpendicular}.1)${}_j$}
\partial^{G}_{j+1-n}\alpha_{j+1}-(-1)^n\alpha_{j}\partial^{F}_{j+1}=0
\quad\text{holds for all}\quad j\in\BZ\,.
 \end{equation}
We need to find a homomorphism $\beta\in\Homb RFG^{n-1}$ that
satisfies
\begin{equation}
\tag*{(\ref{perpendicular}.2)${}_j$}
\partial^{G}_{j+2-n}\beta_{j+1}-(-1)^{n-1}
\beta_{j}\partial^{F}_{j+1}=\alpha_{j+1}
\quad\text{for all}\quad j\in\BZ\,.
 \end{equation}

Set $\beta_j=0$ for $j\ge s$ and assume by descending induction on
$j$ that we already have  maps $\beta_j$ satisfying
(\ref{perpendicular}.2)${}_j$ for some integer $i\in[n,s]$
and all $j\ge i$.
Set $\delta^h=(-1)^{n+1}\Hom{R}{\partial^{F}_{h+1}}{G_{h-n}}$ for
each $h$.
Using (\ref{perpendicular}.1)${}_i$ and
(\ref{perpendicular}.2)${}_i$ we get
 \begin{align*}
\delta^{i}\big(\alpha_i-\partial_{i+1-n}^G\beta_i\big)
&=(-1)^{n+1}\alpha_i\partial^{F}_{i+1}
+(-1)^n\partial_{i+1-n}^G\beta_i\partial^{F}_{i+1}\\
&=-\partial^G_{i+1-n}\alpha_{i+1}
+(-1)^n\partial_{i+1-n}^G \beta_i\partial^{F}_{i+1}\\
&=-\partial_{i+1-n}^G\big(\alpha_{i+1}
+(-1)^{n-1}\beta_i\partial_{i+1}^F\big)\\
&=-\partial_{i+1-n}^G\partial_{i+2-n}^G\beta_{i+1}\\
&=0
 \end{align*}
On the other hand, since one has $\Ext iR{M}{{G}_{i-n}}=0$ the sequence
 \[
\Hom{R}{F_{i-1}}{G_{i-n}}\xra{\ \delta^{i-1}\ }
\Hom{R}{F_{i}}{G_{i-n}}\xra{\ \delta^{i}\ }
\Hom{R}{F_{i+1}}{G_{i-n}}
 \]
is exact,  so there exists a homomorphism $\beta_{i-1}\col
F_{i-1}\to G_{i-n}$, such that
 \[
\alpha_{i}-\partial^{G}_{i+1-n}\beta_{i}=\delta^{i-1}(\beta_{i-1})
=-(-1)^{n-1}\beta_{i-1}\partial^{F}_{i}\,.
 \]
Thus, $\beta_{i-1}$ satisfies (\ref{perpendicular}.2)${}_{i-1}$,
so the induction step is complete.  As a result, for each $j\ge
n-1$ we now have a homomorphism $\beta_j\col F_{j}\to G_{j-(n-1)}$
satisfying the equality (\ref{perpendicular}.2)${}_j$. As
$G_{j-(n-1)}=0$ for $j<n-1$, setting $\beta_j=0$ we extend the
equality to all $j\in\BZ$.  We have proved $\Extb nRMN=0$, as
desired.
 \end{proof}

  \begin{chunk}
   \label{complete-res-def}
A \emph{complete resolution} of $M$ is a morphism of complexes $\nu\col
T\to F$ such that $\nu_i$ is bijective for all $i\gg0$, each $T_i$ is
projective, and for all $n\in\BZ$ and every projective $R$-module $P$
one has $\HH nT=0=\HH n{\Hom RTP}$; see \cite[(1.1)]{CK}. (In some contexts it is
assumed, in addition, that the $R$-modules $T_n$ are also finite; no
such hypothesis is needed or made here.) When such a complete resolution
exists, $\HH{-n}{\Hom RMN}$ is called the $n$th  \emph{Tate cohomology} 
of $M$ with coefficients in $N$.
 \end{chunk}
 
Cornick and Kropholler \cite[(1.2)]{CK} prove that when Tate cohomology is 
defined it is naturally isomorphic to stable cohomology.
We deduce this from Theorem \ref{perpendicular}:

\begin{corollary}
 \label{complete-res}
If $\nu\col T\to F$ is a complete resolution of $M$, then one has
 \[
\CH n{\Hom RTN}\cong\Extv nRMN\quad\text{for each}\quad n\in\BZ\,.
 \]
 \end{corollary}

 \begin{proof}
Fix $n\in\BZ$ and set $K=\Coker(\dd^{T}_{n})$.  For each $i\ge 1$
one then has
 \begin{equation}
\CH{n-1+i}{\Hom RTN}=\Ext iRKN\,,
\tag*{(\ref{complete-res}.1)${}_i$}
 \end{equation}\setcounter{equation}{1}%
because $\shift^{-(n-1)}(T_{\ges n-1})$ is a projective resolution
of $K$.  {}From the condition on $T$ one gets $\Ext
iRKP=0$ for all $i\ge1$ and every projective $R$-module $P$, so
 \begin{equation}
 \label{iso:Tate2}
\Ext1RKN\cong\Extv 1RKN\,.
 \end{equation}
holds by the theorem.  Choose $p\ge n$ with $\nu_i$ bijective for
$i\ge p$. The $R$-module $L=\Ker(\dd^F_{p-1})$ is then isomorphic
to $\Ker(\dd^T_{p-1})$, so there exist exact sequences
\begin{gather*}
0\lra L\lra T_{p-1}\lra\cdots\lra T_n\lra T_{n-1}\lra K\lra0
 \\
0\lra L\lra F_{p-1}\lra\cdots\lra F_1\lra F_{0}\lra M\lra0
 \end{gather*}
In view of \eqref{finite-pd}, the iterated connecting maps
defined by these sequences yield
 \begin{equation}
 \label{iso:Tate3}
\Extv 1RKN\cong\Extv{n-p}RLN\cong\Extv nRMN\,.
 \end{equation}
To finish the proof, concatenate the isomorphisms
(\ref{complete-res}.1)${}_1$, \eqref{iso:Tate2}, and
\eqref{iso:Tate3}.
  \end{proof}

 \section{Additional structures}
 \label{Additional structures} 

In this section we discuss the existence of finer natural
structures on stable cohomology groups, such as rings of
operators or internal gradings.

\begin{proposition}
 \label{stable-finite}
Let $R$ be an algebra over a commutative ring $K$.
 \begin{enumerate}[\rm\quad(1)]
  \item
The exact sequence \eqref{sequence:homology1} is one of graded $K$-modules,
and the various pairings of cohomology groups are $K$-bilinear.
 \end{enumerate}
If, in addition, $K$ is noetherian, $R$ is finite as a
$K$-module, $M$ and $N$ are finite $R$-modules, and $n$ is an integer, then
the following assertions also hold.
 \begin{enumerate}[\rm\quad(1)]
  \item[\rm(2)]
The $K$-modules $\Extv{n}R{M}{N}$ and $\Extb{n+1}R{M}{N}$ are finite
simultaneously.
 \item[\rm(3)]
When $\Ext{\gg j}R{M}{R}=0$ the $K$-module $\Extv nRMN$ is finite
for every $n\in\BZ$.
 \end{enumerate}
 \end{proposition}

 \begin{proof}
(1)  This is due to  the fact that the relevant maps in cohomology are
induced by morphisms of complexes or $K$-modules.

(2)  By (1), the maps in the exact sequence \eqref{sequence:homology2} are
$K$-linear, and under our hypotheses the $K$-modules
$\Ext{n}R{M}{N}$ and $\Ext{n+1}R{M}{N}$ are noetherian.

(3) Let $F\to M$ and $G\to N$ be resolutions by finite projective
$R$-modules.  Choose $m\ge1$ so that 
$\HH{-j}{\Hom RFR}= \Ext{j}R{M}{R}=0$ for all $j>m$; then
$\Hom RFR$ is quasi-isomorphic to the complex $C$ of 
right $R$-modules defined by
 \[
C_i=\begin{cases}
\Hom R{F_{-i}}R&\text{ for }0\ge i\ge -m\,;\\
\im\Hom R{\dd^F_{m+1}}R&\text{ for }i=-m-1\,;\\
0&\text{ for }i<-m-1 \text{ or } i>0\,.
 \end{cases}
 \]
As a consequence, $\Hom RFR\otimes_R G$ is quasi-isomorphic to $C\otimes_R
G$.  For each $n\in\BZ$ this gives the second isomorphism below; Lemma
\ref{omega} provides the first one:
 \begin{align*}
\Extb nRMN&=\CH n{\Homb RFG}\\
&\cong\CH n{{\Hom RFR}\otimes_RG}\\
&\cong\CH n{C\otimes_RG}
 \end{align*}
As $(C\otimes_RG)^n=\coprod_{j=0}^m C^j\otimes_R G^{n-j}$ is a
finite $K$-module for each $n$, we see that $\Extb nRMN$ is a
finite $K$-module; by (2), so is $\Extv nRMN$.
 \end{proof}

Stable cohomology behaves predictably under flat base change.

\begin{proposition}
 \label{flat-extension}
Let $R$ be a commutative noetherian ring, $M$ an $R$-module that 
admits a resolution by finite projective modules, and $R\to R'$ a
homomorphism of rings such that the right $R$-module $R'$ is flat.

For each $R$-module $N$ there is then a commutative diagram
\[
 \xymatrixrowsep{2pc}\xymatrixcolsep{2.5pc}\xymatrix{
R'\otimes_R\Ext{}RMN \ar@{->}[r]^-{R'\otimes_R\iota}\ar@{->}[d]_{\cong}
 &R'\otimes_R\Extv{}RMN \ar@{->}[d]^{\cong}
\\
\Ext{}{R'}{R'\otimes_RM}{R'\otimes_RN}\ar@{->}[r]^-{\iota}
 &\Extv{}{R'}{R'\otimes_RM}{R'\otimes_RN}
 }
\]
When $N=M$ all the maps in the diagram are homomorphisms of 
graded algebras.
 \end{proposition}

 \begin{proof}
Set $(-)'=(R'\otimes_R-)$.  Let $F\to M$ be a resolution by finite
projective $R$-modules and $G\to N$ be a projective resolution.  
In the commutative square
\[
 \xymatrixrowsep{2pc}\xymatrixcolsep{2pc}\xymatrix{
R'\otimes_R\Hom R{F}{G}\ar@{->}[d]\ar@{->}[r]^{\simeq}
&R'\otimes_R\Hom R{F}{N}\ar@{->}[d]^{\cong}
\\
\Hom {R'}{F'}{G'}\ar@{->}[r]\ar@{->}[r]^{\simeq} &\Hom
{R'}{F'}{N'}
 }
 \]
the isomorphism is due to the choice of $F$.  Thus, the left vertical map
is a quasi-isomorphism.  It appears in the following commutative diagram 
where the vertical arrows are induced by the map $\alpha\mapsto\alpha'$
and the rows are exact, see \eqref{sequence:complexes}:
 \[
\xymatrixrowsep{2pc}\xymatrixcolsep{1pc}\xymatrix{ 0\ar@{->}[r]
&R'\otimes_R\Homb R{F}{G}\ar@{->}[r]\ar@{->}[d]_{\cong}
&R'\otimes_R\Hom R{F}{G}\ar@{->}[d]_{\simeq}\ar@{->}[r]
&R'\otimes_R\Homv R{F}{G}\ar@{->}[r]\ar@{->}[d] &0
\\
0\ar@{->}[r] &\Homb {R'}{F'}{G'}\ar@{->}[r] &\Hom
{R'}{F'}{G'}\ar@{->}[r] &\Homv {R'}{F'}{G'}\ar@{->}[r] &0
 }
 \]
The flatness of $R'$ implies that $F'\to M'$ and $G'\to N'$ are 
$R'$-projective resolutions, and that the homology of the right 
hand square above is the desired diagram.
 \end{proof}
\end{subsection}

Next we turn to cohomology of graded objects.

\begin{chunk}
 \label{gradings}
We say that the ring $R$ is \emph{internally graded} if
$R=\bigoplus_{i\in\BZ}^\infty R_i$ as abelian groups, and 
$R_iR_j\subseteq R_{i+j}$ holds for all $i,j$. Internal gradings 
for $M,N$ are defined similarly.  By convention, we allow $M_i$ to 
be written also as $M^{-i}$.  As
usual, we let $M(s)$ denote the graded $R$-module with $M(s)_i=M_{s+i}$
for all $i\in\BZ$.

Assume $R$, $M$, and $N$ are internally graded.  A homomorphism
$\beta\col M\to N$ is \emph{homogeneous} of internal degree $-j$
if $\beta(M_i)\subseteq N_{i-j}$ holds for each $i\in\BZ$. All
such maps form an abelian subgroup $\Homgr RMN^j$ of $\Hom RMN$.
Clearly, the sum of these subgroups is direct, so $\Hom RMN$
contains as a subgroup the group
 \[
\Homgr RMN=\bigoplus_{j\in\BZ}\Homgr RMN^j\,.
 \]
When $M$ is finitely presented, one has $\Homgr RMN=\Hom RMN$.

Every graded $R$-module $M$ has a \emph{graded free resolution}
$F\to M$ that is, a resolution in which each $F_i$ is a graded
free $R$-module and differentials are homogeneous of internal
degree $0$. It produces a subcomplex $\Homgr RFN$ of $\Hom RFN$,
consisting of graded abelian groups and homomorphisms of internal
degree $0$.

Assume that each $R$-module $F_i$ is finite.  One then has 
$\Homgr RFN=\Hom RFN$, and
hence the absolute Ext groups inherit an internal grading:
 \[
\Ext nRMN=\bigoplus_{j\in\BZ}\Ext nRMN^{j}\,.
 \]
 \end{chunk}

For each $n\in\BZ$ one also has equalities
 \[
\Hom RFG_n=\prod_{i\in\BZ}\bigoplus_{j\in\BZ}\Homgr R{F_i}{G_{i+n}}^j\,.
 \]
However, there is no induced internal grading on the right hand side,
so extra steps are needed to introduce such a grading on stable 
cohomology groups.

\begin{proposition}
 \label{stable-graded}
Assume $R$ is an internally graded ring, $M,N$ are internally
graded $R$-modules, and $M$ has a graded free resolution 
$F\to M$ where each $F_i$ is finite (as is the case, for example,
when $R$ is left noetherian and $M$ is finite).

For each $n\in\BZ$ the abelian groups $\Extv{n}R{M}{N}$ and
$\Extb{n}R{M}{N}$ then have natural internal gradings, which are
preserved by the homomorphisms in the exact sequence
\eqref{sequence:homology2} and are additive under the various products.

If $M$ has a complete resolution $T$ by finite projective graded $R$-modules
with differentials  $\dd^T$ of degree $0$, then the internal
gradings of $\CH n{\Homgr RTN}$ and $\Extv nRMN$ are preserved by
the isomorphisms of Corollary \emph{\ref{complete-res}}.
 \end{proposition}
 
 \begin{proof}
Let $G\to N$ be graded free resolution.  For all $j,n\in\BZ$ the
subgroups
 \[
\Homgr RFG^j_n=\prod_{i\in\BZ}\Homgr R{F_i}{G_{i+n}}^j
 \]
of $\Hom RFG_n$ form a subcomplex $\Homgr RFG^j$ of $\Hom RFG$.
The Comparison Theorem for graded resolutions shows that the
canonical morphism $\Homgr RFG^j\to\Homgr RFN^j$ is a
quasi-isomorphism. It follows that the complex $\Homgr
RFG=\bigoplus_{j\in\BZ}\Homgr RFG^j$ appears in a commutative
diagram
\[
 \xymatrixrowsep{2pc}\xymatrixcolsep{2pc}\xymatrix{
\Homgr RFG\ar@{->}[d]\ar@{->}[r]^{\simeq}
 &\Homgr RFN\ar@{=}[d]
 \\
\Hom RFG\ar@{->}[r]^{\simeq} &\Hom RFN
 }
 \]
where the horizontal maps are quasi-isomorphisms, and the equality
is due to the finiteness of the modules $F_i$.  Thus, the left
vertical map is a quasi-isomorphism.

Setting $\Hombgr RFG_n=\bigoplus_{j\in\BZ}\coprod_{i\in\BZ}\Homgr
R{F_i}{G_{i+n}}^j$ for each $n\in\BZ$ one gets an internally
graded subcomplex $\Hombgr RFG$ of the internally graded complex
$\Homgr RFG$.  Thus, $\Homvgr RFG=\Homgr RFG/\Hombgr RFG$ is an
internally graded complex.  On the other hand, one has equalities
 \begin{align*}
\Hombgr RFG_n
 &=\bigoplus_{j\in\BZ}\coprod_{i\in\BZ}\Homgr
R{F_i}{G_{i+n}}^j
 \\
 &=\coprod_{i\in\BZ}\bigoplus_{j\in\BZ}\Homgr
R{F_i}{G_{i+n}}^j
 \\
&=\coprod_{i\in\BZ}\Hom R{F_i}{G_{i+n}}=\Homb RFG_n
 \end{align*}
Putting the preceding observations together, one gets a
commutative diagram
 \[
\xymatrixrowsep{2pc}\xymatrixcolsep{1pc}\xymatrix{
 0\ar@{->}[r]
 &\Hombgr RFG\ar@{->}[r]\ar@{=}[d]
 &\Homgr RFG\ar@{->}[d]_{\simeq}\ar@{->}[r]
 &\Homvgr RFG\ar@{->}[r]\ar@{->}[d]
 &0
\\
 0\ar@{->}[r]
 &\Homb RFG\ar@{->}[r]
 &\Hom RFG\ar@{->}[r]
 &\Homv RFG\ar@{->}[r]
 &0
 }
 \]
It implies that the vertical arrow on the right hand side is a
quasi-isomorphism.

The first assertion of the proposition follows. The remaining ones are
easy consequences of the quasi-isomorphisms above and the
definition of products, respectively, the definition of the
isomorphisms $\CH n{\Hom RTN}\cong\Extv nRMN$.
 \end{proof}

\section{Non-zero-divisors}
\label{Non-zero-divisors}

A staple in basic homological algebra are `change-of-rings theorems' that
track the behavior of cohomology groups under passage to quotients
modulo non-zero-divisors.  An essential ingredient in such results is
the functoriality of absolute Ext groups in the ring argument.  Stable
cohomology does not enjoy a similar property, so we approach 
change of rings through the natural homomorphisms 
 \[
\iota\col \Ext{}RMN\lra\Extv{}RMN
 \]
from \eqref{sequence}.  To simplify notation, we let
$\widehat{\alpha}$ denote the image of $\alpha$ under $\iota$.

For absolute cohomology the following result is due to Gulliksen \cite{Gu}.

\begin{theorem}
 \label{central}
Let $Q$ be an associative ring, $f\in Q$ a central non-zero-divisor,
and set $R=Q/(f)$.  For all $R$-modules $M$, $N$ there exist elements
\[
\vartheta^M\in\Ext2RMM \quad\text{and}\quad
\vartheta^N\in\Ext2RNN
\]
with the following properties:
  \begin{enumerate}[\quad\rm(1)]
 \item
Every $\xi$ in $\Extb{}RMN$ (respectively, $\Ext{}RMN$,
$\Extv{}RMN$) satisfies
\begin{equation*}
\vartheta^N\cdot\xi=\xi\cdot\vartheta^M\,.
\end{equation*}
 \item
$\vartheta^M$ is in the center of $\Ext{}RMM$ and
$\vartheta^N$ in that of $\Ext{}RNN$. 
 \item
${\wh\vartheta}^M$ is in the center of $\Extv{}RMM$ and
${\wh\vartheta}^N$ in that of $\Extv{}RNN$.
 \end{enumerate}
  \end{theorem}

When there is no ambiguity, we let $\vartheta$ denote either one of
$\vartheta^N$ or $\vartheta^M$, and $(\ )_\vartheta$ the functor 
of graded localization at the multiplicatively closed set 
$\{\vartheta^i\var i\ge0\}$.

The last assertion of the corollary is due to Buchweitz, see
\cite[(10.2.3)]{Bu}.

\begin{corollary}
 \label{localization}
There are induced structures of  graded algebras on $\Ext{}RNN_\vartheta$
and $\Ext{}RMM_\vartheta$ (respectively, $\Extv{}RNN_\vartheta$
and $\Extv{}RMM_\vartheta$), and an induced left-right-bimodule
structure on $\Ext{}RMN_\vartheta$ (respectively,
$\Extv{}RMN_\vartheta$).
 
The map $\iota$ from\eqref{sequence} induces isomorphisms
of graded algebras
\begin{align*}
\Ext{}RNN_\vartheta&\lra\Extv{}RNN_\vartheta
\\
\Ext{}RMM_\vartheta&\lra\Extv{}RMM_\vartheta
\end{align*}
and an equivariant isomorphism of graded
bimodules over them,
\[
\Ext{}RMN_\vartheta\lra\Extv{}RMN_\vartheta
\]

If $M$ or $N$ has finite projective dimension over $Q$, then
\[
\Extv{}RMN_\vartheta=\Extv{}RMN\,.
\]
 \end{corollary}

The notation of the theorem is in force for the rest of this section.
In the proofs, presented later in this section, we use a construction 
of Eisenbud \cite[(1.1)]{Ei}:

\begin{chunk}
 \label{eisenbud}
Let $F$ be a \emph{liftable} complex of projective $R$-modules, meaning
that there is a graded projective $Q$-module $\wt F$ with $R\otimes_Q\wt F=F$; 
for example, every complex of free modules is liftable. Since each
$\wt F_i$ is projective, one can choose a map
 \[
\wt\dd^F\in\Hom Q{\wt F}{\wt F}_{-1} \quad\text{with}\quad
R\otimes_Q\wt\dd^F=\dd^F\,.
 \]
As $R\otimes_Q({\wt\dd}^F)^2=(R\otimes_Q\wt\dd^F)^2=0$, for each $x\in\wt
F_n$ there exists $y\in\wt F_{n-2}$ satisfying $({\wt\dd}^F)^2(x)=fy$. As
$f$ is a non-zero-divisor on $\wt F_{n-2}$ (it is one on $Q$ and
$\wt F_{n-2}$ is projective), $y$ is defined uniquely and hence depends
$Q$-linearly on $x$. Setting $\wt\theta^F(x)=y$ one gets a homomorphism
$\wt\theta^F\in\Hom Q{\wt F}{\wt F}_{-2}$.  As $f$ is central, we get
 \[
f(\wt\theta^F\wt\dd^F)=(f\wt\theta^F)\wt\dd^F
 =({\wt\dd}^F)^3=\wt\dd^F(f\wt\theta^F)
=f(\wt\dd^F\wt\theta^F)\,.
 \]
As $f$ is a non-zero-divisor, this implies
$\wt\theta^F\wt\dd^F=\wt\dd^F\wt\theta^F$.  Thus, one gets a chain map
 \[
\theta^F=R\otimes_Q\wt\theta^F\in\Hom R{F}{F}_{-2}\,.
 \]
  \end{chunk}

The first assertion of the next result is \cite[(1.3)]{Ei}.  An adaptation of the 
original argument allows us to handle the other two cases as well.

\begin{lemma}
\label{homotopy}
Let $F$ and $G$ be liftable complexes of projective $R$-modules and let 
$\gamma$ be a homomorphism in $\Hom RFG$.

If $\gamma$ is a chain map, then $\theta^G\gamma$ and $\gamma\theta^F$
are homotopic in $\Hom RFG$.

If $\gamma$ is a bounded chain map, then $\theta^G\gamma$ and
$\gamma\theta^F$ are homotopic in $\Homb RFG$.

If $\wh\gamma$ is a chain map, then $\theta^G\wh\gamma$ and
$\wh\gamma\theta^F$ are homotopic in $\Homv RFG$.
 \end{lemma}

\begin{proof}
Set $n=\lfloor \gamma\rfloor$ and assume first $\wh\gamma$ is a chain
map. Thus, the map
 \begin{equation}
\label{almost}
\delta=\dd^G\gamma-(-1)^{n}\gamma\dd^F+\delta\in\Hom RFG_{n-1}
  \end{equation}
satisfies $\delta_i=0$ for all $i\ge j$
and a fixed $j\in\BZ$.  Choose $\wt\gamma\in\Hom Q{\wt F}{\wt G}_n$ with
$R\otimes_Q\wt\gamma=\gamma$ and $\wt\delta\in\Homb Q{\wt F}{\wt G}_{n-1}$
with $R\otimes_Q\wt\delta=\delta$ and $\wt\delta_i=0$ for all $i\ge j$.
There exits then a unique $\wt\tau\in\Hom Q{\wt F}{\wt G}_{n-1}$ satisfying
 \[ 
\wt\dd^G\wt\gamma=(-1)^{n}\wt\gamma\wt\dd^F+\wt\delta +f\wt\tau\,.
  \]
Using the relation above and the equalities $({\wt\dd}^F)^2=f\wt\theta^F$
and $({\wt\dd}^G)^2=f\wt\theta^G$, we get
 \begin{align*}
f(\wt\theta^G\wt\gamma)
&=({\wt\dd}^G)^2\wt\gamma\\
&=(-1)^{n}\wt\dd^G\wt\gamma\wt\dd^F+\wt\dd^G\wt\delta+\wt\dd^G f\wt\tau\\
&=\wt\gamma({\wt\dd}^F)^2+(-1)^{n}\wt\delta\wt\dd^F
 +(-1)^{n}f\wt\tau\wt\dd^F+\wt\dd^G\wt\delta+f\wt\dd^G\wt\tau\\
&=f\big(\wt\gamma\wt\theta^F+\wt\dd^G\wt\tau+(-1)^{n}\wt\tau\wt\dd^F\big)
 +\big(\wt\dd^G\wt\delta+(-1)^{n}\wt\delta\wt\dd^F\big)\,.
  \end{align*}
Since $f$ is a non-zero-divisor on $\wt G$, the preceding computation
yields
 \[
\wt\theta^G\wt\gamma-\wt\gamma\wt\theta^F-
\big(\wt\dd^G\wt\tau-(-1)^{n+1}\wt\tau\wt\dd^F\big)\in\Homb Q{\wt
F}{\wt G}\,.
 \]
The map $\tau=R\otimes_Q\wt\tau\col F\to G$ then satisfies
 \[
\theta^G\gamma-\gamma\theta^F-
\big(\dd^G\tau-(-1)^{n+1}\tau\dd^F\big)\in\Homb RFG\,.
 \]
In other words, $\wh\tau\in\Homv RFG$ is a homotopy between
$\theta^G\wh\gamma$ and $\wh\gamma\theta^F$.

If $\gamma$ is a chain map, then \eqref{almost} holds with $\delta=0$,
so in the computation above one can choose $\wt\delta=0$.  The resulting
$\tau$ is a homotopy between $\theta^G\gamma$ and $\gamma\theta^F$.

When $\gamma$ is a bounded chain map the map $\wt\gamma$ can be 
chosen to be bounded as well, and then $\wt\tau$ is necessarily bounded, 
so the homotopy $\tau$ is in $\Homb RFG$.
 \end{proof}

\begin{proof}[Proof of Theorem \emph{\ref{central}}] 
Let $F\to M$ be a liftable projective resolution (for example, 
choose $F$ to be free a resolution), and set $\vartheta^M=
\cls{\theta^F}\in\Ext2RMM$.

First we show that $\vartheta^M$ does not depend on the choice of 
liftable resolution.  If $F'\to M$ is one, then pick
a morphisms of complexes $\gamma\col F\to F'$ lifting $\id M$.  By 
Lemma \ref{homotopy}, the maps $\theta^{F'}\gamma$ and  
$\gamma\theta^F$ are homotopic.  Thus, the isomorphisms
 \[
\HH{}{\Hom R{F'}{F'}}\xra{\ \cong\ }\HH{}{\Hom R{F}{F'}}
\xla{\ \cong\ }\HH{}{\Hom R{F}{F}}\,.
 \]
map $\cls{\theta^{F'}}$ and $\cls{\theta^F}$  to the same element,
which was to be shown.

The other assertions of the theorem follow directly from the lemma.
 \end{proof}

One can produce liftable resolutions using a construction of Shamash 
\cite[\S 3]{Sh}; we describe it next, following the simplified exposition in \cite{res}.

\begin{chunk}
 \label{shamash}
Let $E$ be a projective resolution of $M$ over $Q$.

By induction, one gets for each $i\ge0$ a map $\sigma^{(i)}\in \Hom
QEE_{2i-1}$, such that
 \[
\sigma^{(0)}=\dd^E
 \quad\text{and}\quad
\sum_{h=0}^i\sigma^{(h)}\sigma^{(i-h)}=
 \begin{cases}
f\id E&\text{ for}\quad i=1\,;\\0&\text{ for}\quad i\ge2\,.
 \end{cases} \]

Let $D$ be a graded $\BZ$-module, such that for each $i\ge0$ the
$\BZ$-module $D_{2i}$ is free with a single basis element
$y^{(i)}$.  Let $\wt F$ be the graded projective $Q$-module with
 \[
\wt F_n=\bigoplus_{i\ges0}E_{n-2i}\otimes_\BZ D_{2i}\,.
 \]
For every $i\ge0$ and each $e\in E_{n-2i}$ the formula
\[
\wt\dd_n(e\otimes y^{(i)})= \sum_{h\ges0}\sigma^{(h)}(e)\otimes
y^{(i-h)}
\]
defines a $Q$-linear map $\wt\dd\col \wt F\to\wt F$ of degree
$-1$. A direct computation yields
\[
\wt\dd{}^2(e\otimes y^{(i)})=fe\otimes y^{(i-1)}\,.
\]
As a consequence, one obtains a complex of projective $R$-modules
 \[
(F,\dd)=(R\otimes_Q\wt F,R\otimes_Q\wt\dd)\,.
 \]

It is proved in \cite[(3.1.3)]{res} that this a resolution of $M$ over $R$.
This  projective resolution is clearly liftable, and one can define a map 
$\theta^F$ as in \eqref{eisenbud} by setting 
 \[
\theta^F(r\otimes e\otimes y^{(i)})=r\otimes e\otimes y^{(i-1)}\,.
 \] 
 \end{chunk}
 
\begin{proof}[Proof of Corollary \emph{\eqref{localization}}]
The multiplicativity properties follow from the definitions of the various
products in cohomology and the centrality of the element $\vartheta$.

The exact sequence \eqref{sequence:homology1} is the homology 
sequence of an exact sequence of DG modules over $\Hom RNN$.
Thus, its maps commute with left multiplication by the element
$\vartheta=\vartheta^N\in\Ext2RNN$.  Localizing at $\{\vartheta^i\}_{i\ges0}$
the exact sequence \eqref{sequence} we obtain an exact sequence of 
$\Ext{}RNN{}_\vartheta$-$\Ext{}RMM{}_\vartheta$ bimodules.  
Lemma \ref{nilpotent} implies that  in this sequence $\Extb{}RMN_\vartheta$ 
vanishes,  so $\iota_\vartheta$ is bijective.

For the last assertion, we show that if $\pd QM$ is finite, then
$\wh\vartheta\in\Extv{}RMM$ is invertible.  Let $E$ be a
finite projective resolution of $M$ over $Q$ and set
\[
\zeta(r\otimes e\otimes y^{(i)})=r\otimes e\otimes y^{(i+1)}
\]
in the notation of \eqref{shamash}. This formula defines a map
$\zeta\in\Hom QFF_2$, such that
\[
(\dd\zeta-\zeta\dd)(F)\subseteq\big(R\otimes_QE\otimes_\BZ \BZ
y^{(0)}\big) \supseteq(\zeta\theta^F-\id F)(F)\quad\text{and}\quad
\theta^F\zeta=\id F\,.
\]
Thus, $\dd(\wh\zeta)=0$ and $\wh\zeta\wh\theta^F=\wh\theta^F\wh\zeta=
\wh{\id F}$ in $\Homv RFF$, as desired.
 \end{proof}

As another application of Shamash's construction, we derive an
exact sequence used several times in the paper, for which many
other proofs are known.

\begin{proposition}
 \label{change-of-rings}
There is an exact sequence of graded bimodules
 \[
 \xymatrixrowsep{.5pc} \xymatrixcolsep{2pc}\xymatrix{
\shift^{-2}\Ext{}{R}MN \ar@{->}[r]^-{\lambda}&\Ext{}{R}MN
\ar@{->}[r]&\Ext{}QMN\ar@{->}[r]&
\\
\shift^{-1}\Ext{}{R}MN\ar@{->}[r]^-{\lambda}&\shift\Ext{}{R}MN
 }
 \]
over $\Ext{}RNN$-$\Ext{}RMM$, where $\lambda$ is given by
multiplication with $\vartheta$.
 \end{proposition}

\begin{proof}
Using the notation of \eqref{shamash}, we form a sequence
\[
0\lra R\otimes_QE \xra{\ \alpha\ }F\xra{\ \theta^F\
}\shift^2F\lra0
\]
of morphisms complexes of $R$-modules, where $\alpha(r\otimes e)=
r\otimes e\otimes y^{(0)}$.  It is clear that the underlying
sequence of graded $R$-modules is split exact.  Let $G\to N$ be a
projective resolution.  The induced sequence of complexes of
abelian groups
\[
0\lra\shift^{-2}\Hom RFG\lra\Hom RFG\lra\Hom R{R\otimes_QE}G\lra0
\]
is then exact, and its cohomology exact sequence is the desired
one.
 \end{proof}

\section{Depth of cohomology modules}
\label{Depth of cohomology modules}

Let $(R,\fm,k)$ be a local ring\footnote{Recall that this
means that $R$ is a commutative noetherian ring with unique maximal 
ideal $\fm$ and residue field $k=R/\fm$. },  $M$ a finite $R$-module, and set
 \[
\CE=\Ext{}{R}kk\quad\text{and}\quad\CM=\Ext{}RMk\,.
\]
The \emph{depth} of $\CM$ over $\CE$ is defined by means of the
formula
\[
\depth{\CE}\CM=\inf\{n\in\BN\mid\Ext n\CE k\CM\ne0\}\,.
\]
We systematically write $\depth{}\CE$ in place of
$\depth{\CE}\CE$.

The use of depth to study the structure of $\CE$ was pioneered by
F\'elix, Jacobsson, Halperin, L\"ofwall, and Thomas in the important
paper \cite{FHJLT}.  Their main result is the finiteness of
$\depth{}{\CE}$; to prove it they develop methods for obtaining
upper bounds on depth.  To study stable cohomology we
mostly need lower bounds.  

General properties of depth of graded modules used in the paper are 
collected in Appendix \ref{Depth over graded algebras}.  In this section 
we focus on additional properties stemming from the cohomological 
nature of $\CE$ and $\CM$.

\begin{subsection}{Universal enveloping algebras}
 \label{universal}
Let $\pi$ be a graded Lie algebra over $k$, such that
$\rank_k\pi^i$ is finite for all $i\in\BZ$ and $\pi^i=0$ for
$i\le0$, and let $\cD$ be the universal enveloping algebra of
$\pi$; for definitions of these notions see
\cite[(10.1.2)]{res}\footnote{Where condition (3) contains a typo:
both $-$ signs should be changed to $+$ signs.}.

\begin{subchunk}
 \label{PBW}
The $k$-algebra $\cD$ has an increasing multiplicative filtration,
whose $p$th stage is the $k$-linear span of products involving at
most $p$ elements of $\pi$. By the Poincar\'e-Birkhoff-Witt
Theorem, see \cite[Thm.\ 2, Cor.]{Sj2}, the associated graded $k$-algebra is
isomorphic to ${\mathsf\Lambda}_k(\pi^{\operatorname{odd}})
\otimes_k {\mathsf S}_k(\pi^{\operatorname{even}})$, where
${\mathsf\Lambda}_k$ and $\mathsf S_k$ denote, respectively,
exterior algebra and symmetric algebra functors over $k$.
 \end{subchunk}

Directly from the Poincar\'e-Birkhoff-Witt isomorphism one gets:

\begin{subchunk}
 \label{lie-series}
There is an equality of formal power series
 \[
\sum_{n=0}^\infty(\rank_k\cD^n)t^n =
  \frac{\prod_{i\ges0}(1+t^{2i+1})^{\rank_k\pi^{2i+1}}}
   {\prod_{i\ges0}(1-t^{2i+2})^{\rank_k\pi^{2i+2}}}
 \]
 \end{subchunk}

\begin{subchunk}
 \label{PBW-nzd}
If $\cD'$ is the universal enveloping algebra of a graded Lie 
subalgebra $\pi'$ of $\pi$, then $\cD$ is free as a left 
$\cD'$-module and as a right $\cD'$-module.  In particular, every
$\zeta\in\pi^{\operatorname{even}}\smallsetminus\{0\}$ is a
left non-zero-divisor and a right non-zero-divisor.
 \end{subchunk}

A theorem of Milnor and Moore, Andr\'e, and Sj\"odin, see
\cite{Sj2} or \cite[(10.2.1.5)]{res} for proofs, introduces graded
Lie algebras into the study of local rings:

\begin{subchunk}
 \label{lie}
The graded $k$-algebra $\CE=\Ext{}Rkk$ is the universal enveloping algebra
of a graded Lie algebra, denoted $\pi_R$ and called the \emph{homotopy
Lie algebra} of $R$.
 \end{subchunk}

In low degrees the components of $\pi_R$ are easy to describe.

\begin{subchunk}
 \label{cohen}
A \emph{minimal Cohen presentation} of the $\fm$-adic completion $\wh
R$ of $R$ is an isomorphism $\wh R\cong Q/\fa$, where $(Q,\fn,k)$ is a
complete regular local ring and $\fa$ is an ideal contained in $\fn^2$.
Cohen's Structure Theorem shows that one always exists.

There are isomorphisms of $k$-vector spaces, see e.g.\
\cite[(10.2.1.2), (7.1.5)]{res}:
 \begin{align}
 \label{cohen1} \tag{\ref{cohen}.1}
\pi^1_R&\cong\Hom k{\fm/\fm^2}k\cong\Hom k{\fn/\fn^2}k\,.
 \\
 \label{cohen2} \tag{\ref{cohen}.2}
\pi^2_R&\cong\Hom k{\fa/\fn\fa}k\,.
 \end{align}

Recall that the number $\codim R=\edim R-\dim R$ is called the
\emph{codimension}  of $R$.  Krull's Principal Ideal
Theorem and the catenarity of the regular ring $Q$ give
 \begin{equation}
 \label{codimension}\tag{\ref{cohen}.3}
 \rank_k\pi^2_R\ge\operatorname{height}(\fa)=\dim Q-\dim\wh R
 =\edim R-\dim R=\codim R\,.
 \end{equation}
 \end{subchunk}

A specific subalgebra of $\CE$ will prove useful in computations.

\begin{sublemma}
 \label{subalgebra}
The graded subspace $\pi_R^{\ges2}$ of $\pi_R$ is a Lie ideal,
and the universal enveloping algebra $\cD$ of $\pi_R^{\ges2}$
satisfies $\depth{}\cD=\depth{}\CE$.
 \end{sublemma}

 \begin{proof}
For degree reasons, $\pi_R^{\ges2}$ is a graded Lie ideal
of $\pi_R$.  Thus, one has $\cD^{\ges1}\CE=\CE\cD^{\ges1}$.
The Poincar\'e-Birkhoff-Witt Theorem \ref{PBW} implies
that the right
$\cD$-module $\CE$ is free and $\CE/(\CE\cdot\cD^{\ges 1})\cong
{\mathsf\Lambda}_k(\pi^1)$ holds. Now apply Corollary
\ref{extension2}.
  \end{proof}

The following known simple consequence of \eqref{lie} is used throughout
the paper.

\begin{sublemma}
 \label{singular}
The ring $R$ is singular if and only if some element
$\vartheta\in\CE^2$ is a left and right non-zero-divisor on $\CE$,
if and only if $\depth{}\CE\ge1$.
 \end{sublemma}

 \begin{proof}
When $R$ is singular so is $\wh R$, hence $\pi_R^2\ne0$ by
\eqref{cohen2}, so \eqref{PBW-nzd} yields a left and right
non-zero-divisor $\vartheta\in\CE^2$; in view of \eqref{depth0}
the existence of $\vartheta$ implies $\depth{}\CE\ge1$.  On the
other hand, $\depth{}\CE\ge1$ implies $R$ is singular: assuming
the contrary one gets $\CE^i=0$ for $i\gg0$, hence
$\gam\CE=\CE\ne0$, contradicting \eqref{depth0}.
 \end{proof}

The next result is due to Martsinkovsky \cite[Thm.\ 6]{Ma2}. 
We give a short proof. 

\begin{subtheorem}
 \label{martsinkovsky}
When $R$ is singular for each $n\in\BZ$ there is an exact sequence 
 \[
0\lra{\Ext{n}Rkk}\xra{\ \iota\ }{\Extv{n}Rkk}\xra{\ \eth\ }
{\Extb{n-1}Rkk}\lra0
 \] 
 \end{subtheorem}
 
  \begin{proof}
We prove $\eta=0$ in the exact sequence of graded
left $\Ext{}Rkk$-modules
 \[
{\Extb{}Rkk}\xra{\ \eta\ }{\Ext{}Rkk}\xra{\ \iota\ }{\Extv{}Rkk}\xra{\ \eth\ }
\shift{\Extb{}Rkk}\xra{\ \shift\eta\ }\shift{\Extv{}Rkk}
 \]
given by \eqref{sequence}.  In view of \eqref{depth0}, Lemmas 
\ref{nilpotent} and \ref{singular} imply the equalities 
 \[
{\Extb{}Rkk}=\gam{\Extb{}Rkk}
 \quad\text{and}\quad
\gam{\Ext{}Rkk}=0\,,
 \]
where $\gam{}$ denotes the section functor \ref{gamma}. The 
$\Ext{}Rkk$-linearity of $\eta$ now yields
  \[
\eta({\Extb{}Rkk})=\eta(\gam{{\Extb{}Rkk}})\subseteq\gam{\Ext{}Rkk}=0\,.
\qedhere
  \]
   \end{proof}
\end{subsection}

 \begin{subsection}{Regular elements}
Fix an element $g\in\fm$, let $M$ be an $R$-module annihilated by
$g$.  We view $M$ also as a module over $R'=R/(g)$, and set:
 \[
\CE'=\Ext{}{R'}kk
\quad\text{and}\quad
\CM'=\Ext{}{R'}Mk\,.
 \]
The canonical homomorphism of rings $R\to R'$ induces a homomorphism 
 \[
\rho^k\col\CE'\to\CE
 \]
of graded $k$-algebras and a $\rho^k$-equivariant homomorphism 
of graded modules
 \[
\rho^M\col\CM'\to\CM\,.
 \]

\begin{subchunk}
 \label{depth-change}
If $g$ is $R$-regular Theorem \ref{central} and Proposition 
\ref{change-of-rings} yield an exact sequence
 \begin{equation*}
\shift^{-2}\CM'\xra{\ \lambda^M\ }\CM'\xra{\ \rho^M\ }
\CM\lra\shift^{-1}\CM'\xra{\ \shift\lambda^M\ }\shift\CM'
 \end{equation*}
of graded left $\CE'$-modules, with $\lambda^M$ induced by multiplication
on $\CM'$ with a central element $\vartheta'\in\CE'^2$.  Thus, one then has
$\Ker(\rho^k)=\vartheta'\CE=\CE'\vartheta'$ and $\Ker(\rho^M)=\vartheta'\CM'$.
 \end{subchunk}

It is useful to know how depth changes when passing from the
$\CE$-module $\CM$ to the $\CE'$-module $\CM'$.  In two cases we
provide complete---and contrasting---answers.

\begin{subproposition}
 \label{hypersection1}
If $g\notin\fm^2$ is $R$-regular, then $\rho^k$ and $\rho^M$ are
injective,
 \[
\depth{}{\CE'}=\depth{}{\CE} \quad\text{and}\quad
\depth{\CE'}{\CM'}=\depth{\CE'}{\CM}\,.
 \]
 \end{subproposition}

 \begin{proof}
It is proved in \cite[(2.8)]{vpd} that $\rho^k$ and $\rho^M$ are
injective. Set $\cD=\rho^k(\CE')$.

The exact sequence in \eqref{depth-change} and \eqref{linked}
yield $\depth{\cD}{\CM'}=\depth{\cD}{\CM}$.

As $\rho^k$ is injective, applied to $M=k$ the same sequence yields
 \[
\sum_{n=0}^\infty(\rank_k\CE^n)t^n=
(1+t)\sum_{n=0}^\infty(\rank_k\cD^n)t^n\,.
 \]
As $\cD$ is the universal enveloping algebra of the Lie subalgebra
$\pi'=\rho^k(\pi_{R'})\subseteq\pi_R$, from \eqref{lie-series} we deduce
$\rank_k\pi^1_R= \rank_k\pi'^1+1$ and $\rank_k\pi^n_R=\rank_k\pi'^n$
for all $n\ge2$.  Thus, $\pi_R=\pi'\oplus k\varepsilon$ for some
$\varepsilon\in\pi^1_R$, so $\pi'$ is an ideal for degree reasons, and
hence $\cD$ is normal in $\CE$.  The sequence in \eqref{depth-change}
also produces an isomorphism $\CE\cong\cD\oplus\shift^{-1}\cD$
of left $\cD$-modules, so Corollary \ref{extension2} gives
$\depth{}{\cD}=\depth{}{\CE}$.
 \end{proof}

\begin{subproposition}
 \label{hypersection2}
If $g\in\fm^2$ is $R$-regular, then the central element $\vartheta'\in
{\CE'}^2$ from \eqref{depth-change} is central also in $\CS$, is 
regular on $\CE'$, and one has
 \[
\CE'/(\vartheta')\cong\CE
 \quad\text{and}\quad
\depth{}{\CE'}=\depth{}{\CE}+1\,.
 \]
When $g\in\fm\ann_RM$ the element $\vartheta'$ is a regular also
on $\CM'$ and one has
 \[
\CM'/\vartheta'\CM'\cong\CM
 \quad\text{and}\quad
\depth{\CE'}{\CM'}=\depth{\CE}{\CM}+1\,.
 \]
 \end{subproposition}

  \begin{proof}
By Theorem \ref{central} the element $\vartheta'$ is central in $\CS$,
and by \cite[(2.8)]{vpd} it is regular on $\CE'$ and $\CM'$.  
The isomorphisms come from \eqref{depth-change}, the
equalities from \eqref{rees}(2). 
 \end{proof}

A last variation on the preceding theme is proved by elementary
arguments:

\begin{sublemma}
If $M=L/gL$ for some $R$-module $L$ and an $L$-regular element
$g\in\fm$, and $\CL$ denotes the $\CE$-module $\Ext{}R{L}k$, then
one has $\depth{\CE}{\CM}=\depth{\CE}{\CL}$.

When $g$ is also $R$-regular, one has $\depth{\CE'}{\CM'}=\depth{\CE'}{\CL}$
as well.
  \end{sublemma}

\begin{proof}
For the first equality apply \eqref{linked} to the exact sequence of graded 
$\CE$-modules
 \[
0\lra\shift^{-1}\CL\lra\CM\lra\CL\lra0
 \]
It is induced by exact sequence $0\to L\xra{g}L \to M\to0$ 
of $R$-modules because one has $\Ext{}R{g\id M}k=0$.  For the 
second equality, note the isomorphism $\CM'\cong\CL$.
 \end{proof}
  \end{subsection}

\begin{subsection}{Residue field extensions}
For depth computations one can sometimes adjust the ring $R$ while
preserving essential homological properties.

\begin{sublemma}
 \label{residue-extension}
When $R\to (R',\fm',k')$ is a flat homomorphism of local rings
with $\fm R'=\fm'$ there is a commutative diagram of homomorphisms
of graded $k'$-algebras
\[
\xymatrixrowsep{2pc}\xymatrixcolsep{2.5pc}\xymatrix{
k'\otimes_k\CE\ar@{->}[r]^-{k'\otimes_k\iota}\ar@{->}[d]_{\cong}
&k'\otimes_k\CS \ar@{->}[d]^{\cong}
\\
\CE'\ar@{->}[r]^-{\iota}&\CS'
 }
\]
where $\CS=\Extv{}{R}{k}{k}$, $\CE'=\Ext{}{R'}{k'}{k'}$ and
$\CS'=\Extv{}{R'}{k'}{k'}$; in particular,
 \[
\depth{}{\CE'}=\depth{}{\CE}
 \quad\text{and}\quad
\depth{\CE'}{\CS'}=\depth{\CE}{\CS}\,.
 \]
 \end{sublemma}

\begin{proof}
Proposition \ref{flat-extension} provides the commutative
diagram.  It implies the equalities of depths, in view of
standard change of rings formulas.
 \end{proof}

Recall that $\edim R$ denotes the \emph{embedding dimension} of $R$, 
that is, the minimal number of generators of $\fm$, and $\mult R$ denotes the 
multiplicity of $R$.

\begin{sublemma}
 \label{weak}
There exists a complete local ring $R'$ with algebraically closed
residue field, $\mult R =\mult R'$, $\edim R'=\edim R-\depth{}R$,
$\depth{}R'=0$, and
 \[ 
\depth{}{(\Ext{}{R'}{k'}{k'})}=\depth{}\CE\,.
 \] \end{sublemma}

\begin{proof}
Let $k'$ be an algebraic closure of $k$.  There always is a flat local
homomorphism $R\to S$ with $S/\fm S=k'$, see \cite[Ch.\ IX, App.,
Thm.\ 2, Cor.]{Bo}. One then has $\depth{}{\Ext{}
{\wh S}{k'}{k'}}= \depth{}\CE$, see Lemma \ref{residue-extension}, and
$\depth{}{\wh S}=\depth{}R$.

If $\depth{}R$ is positive, then $\mult T=\mult S$ for $T=S/(g)$
and some regular element $g\notin(\fn{}\wh S)^2$.  In that case
$\edim T=\edim R-1$ and  $\depth{}T=\depth{}R-1$ also hold, 
while Proposition \ref{hypersection1}
gives $\depth{}{\Ext{}T{k'}{k'}}= \depth{}{\Ext{}{\wh S}{k'}{k'}}$.
 \end{proof}
   \end{subsection}

\section{Finiteness of stable cohomology}
\label{Finiteness of stable cohomology}

In this section $(R,\fm,k)$ is a local ring.  Classical results
characterize ring-theoretical properties of $R$ in terms of
vanishing of absolute Ext modules.  Here we establish analogs for
stable Ext modules.  Remarkably, the key turns out to be a better 
understanding of bounded cohomology.

\begin{theorem}
 \label{bass-numbers}
For each $R$-module $N$ there is an isomorphism
 \begin{equation}
 \label{bass-numbers.1}
\Extb{}RkN\cong\Ext{}RkR\otimes_k\Tor{}RkN
 \end{equation}
of graded left $\Ext{}RNN$-modules.  In particular, there is an
isomorphism
\begin{equation}
 \label{bass-numbers.2}
\Extb{n}RkN\cong\coprod_{i-j=n}^\infty\Ext iRkR\otimes_k\Tor
jRkN
 \end{equation}
of $k$-vector spaces for every $n\in\BZ$.
 \end{theorem}

\begin{noremark}
The reader will notice that the graded algebra $\Ext{}Rkk$ acts from
the right on both $\Extb{}RkN$ and $\Ext{}RkR$, and from the left
on $\Tor{}RkN$. We do not know whether these structures are
related.
 \end{noremark}

\begin{proof} Let $G\to N$ and $F\to k$ be free 
resolutions, with $F_n$ finite for each $n$, and let $R\to J$ be an 
injective resolution. These resolutions induce quasi-isomorphisms
 \begin{align*}
\Hom RFR\otimes_RG\simeq&\Hom RFJ\otimes_RG\\
\simeq&\Hom RkJ\otimes_RG\\
\cong&\Hom RkJ\otimes_k(k\otimes_RG)
 \end{align*}
that commute with the action of $\Hom RGG$. {}From Lemma
\ref{omega} and the K\"unneth formula one now obtains
isomorphisms of graded left $\Ext{}RNN$-modules
  \[
\Extb{}RkN\cong\HH{}{\Hom RFR\otimes_RG}
\cong\Ext{}RkR\otimes_k\Tor{}RkN\,. 
 \qedhere
   \]
   \end{proof}

The \emph{Bass numbers} of $R$, defined by $\mu^n=\rank_k\Ext nRkR$,
appear in some of the most useful characterizations of the Gorenstein
property, recalled below.

 \begin{chunk}
 \label{gorenstein-defs}
The ring $R$ is \emph{Gorenstein} if it satisfies the following equivalent 
conditions:
 \begin{enumerate}[\rm\quad(i)]
 \item
$\mu^n=1$ for $n=\depth{}R$ and $\mu^n=0$ for all $n\ne \depth{}R$.
\item
$\mu^n=0$ for some $n>\depth{}R$.
  \item
$\idim RR<\infty$.
 \end{enumerate}
  \end{chunk}
 
\begin{corollary}
 \label{stable-pd}
If the $k$-vector space $\Extv nRkN$ is finite for some $n\in\BZ$, then
$N$ has finite projective dimension or $R$ is Gorenstein.
 \end{corollary}

\begin{noremark}
It is tempting to ask whether finiteness conditions on $\Extv nRMk$
imply homological restrictions on the module $M$ or the ring $R$.
We have no answer to the first question.  A negative answer to the
second one is given in Example \ref{arbitrary}.
 \end{noremark}

\begin{proof}
When $\rank_k\Extv nRkN$ is finite so is $\rank_k\Extb{n+1}RkN$;
see Proposition \ref{stable-finite}(2).  In view of equation
\eqref{bass-numbers.2} this implies $\Tor jRkN=0$ for all $j\gg0$ or
$\mu^i=0$ for all $i\gg0$; that is, $N$ has finite projective dimension
or $R$ is Gorenstein.
 \end{proof}

Even the entire sequence $(\rank_k\Ext nRkk)_{n\ges0}$ does not recognize
the Gorenstein property of $R$, see Example \ref{ambiguous}, so the next
result is rather surprising.
The expression for $\rank_k\Extv{n}Rkk$ in the next result is known,
see \cite[(9.2)]{AM}.  

\begin{corollary}
 \label{gorenstein}
The ring $R$ is Gorenstein if $\rank_k\Extv{n}Rkk<\infty$ for some
$n\in\BZ$, only if the module $\Extv{n}RMN$ is finite for every 
$n\in\BZ$ and all finite modules $M,N$.

If $R$ is Gorenstein and singular, then for $d=\depth{}R$ and each
$n\in\BZ$ one has
 \[
\rank_k\Extv{n}Rkk=\rank_k\Ext{n}Rkk+\rank_k\Ext{d-1-n}Rkk\,.
\]
  \end{corollary}

 \begin{proof}
If $\rank_k\Extv{n}Rkk$ is finite for some $n\in\BZ$, then $R$ is Gorenstein
by Corollary \ref{stable-pd}.  If $R$ is Gorenstein, then it has finite
injective dimension as $R$-module, see \eqref{gorenstein-defs}, so
by Theorem \ref{stable-finite}(3) the module  $\Extv nRMN$ is finite when
$M$ and $N$ are.

Theorem \ref{martsinkovsky} provides the first equality in the 
following chain:
 \begin{align*}
\rank_k\Extv{n}Rkk
&=\rank_k\Ext{n}Rkk+\rank_k\Extb{n-1}Rkk\\
&=\rank_k\Ext{n}Rkk+\rank_k\Tor{d-1-n}Rkk\\
&=\rank_k\Ext{n}Rkk+\rank_k\Ext{d-1-n}Rkk
 \end{align*}
For the others use formulas \eqref{bass-numbers.2} and 
\eqref{gorenstein-defs}(i), and vector space duality.
 \end{proof}

\begin{corollary}
 \label{regular}
The ring $R$ is regular if $\Extv{n}Rkk=0$ for some $n\in\BZ$,
only if $\Extv nRMN=0$ for every $n\in\BZ$ and all $R$-modules
$M,N$.
  \end{corollary}

\begin{proof}
When $R$ is regular it has finite global dimension, so $\Extv
nRMN=0$ for all $M$, $N$, and $n$ by \eqref{finite-pd}. When $R$
is singular $\Ext{i}Rkk\ne0$ for all $i\ge0$.  Assuming $\Extv{n}Rkk=0$ 
for some $n$, part (3) of the theorem shows that $R$ is Gorenstein.  
Part (2) yields $\Ext{n}Rkk=0=\Ext{d-1-n}Rkk$, implying $d<0$, 
which is absurd.
 \end{proof}

It is well known that polynomial growth of the sequence $(\rank_k\Ext
nRkk)_{n\ges0}$ recognizes a smaller class of rings, whose definition
we proceed to recall.

 \begin{chunk}
  \label{ci-defs}
The ring $R$ is \emph{complete intersection} if in some minimal 
Cohen presentation of $\wh R\cong Q/\fa$, see \eqref{cohen}, the
ideal $\fa$ is generated by a $Q$-regular set.
 \end{chunk}
 
Next we strengthen the recognition criterion for complete intersections.

\begin{theorem}
 \label{ci-betti}
Let $R$ be a local ring, set $d=\dim R$ and $c=\codim R$.

For each $n\in\BZ$ there is an inequality
\begin{equation*}
\label{absolute_betti} \tag*{\rm(\ref{ci-betti}.1)${}_n$}
\rank_k\Ext nRkk\ge\sum_{i=0}^{d}\binom di\binom{c+n-i-1}{c-1}\,.
 \end{equation*}

When $R$ is complete intersection equalities hold for all $n\ge0$.

If equality holds for a single $n\ge2$, then $R$ is complete
intersection.
 \end{theorem}

\setcounter{equation}{1}
\begin{Remark}
 \label{ci-remark1}
For every ring $R$ easy computations yield $\rank_k\Ext0Rkk=1$ and
$\rank_k\Ext1Rkk=\edim R=c+d$. This shows that in (\ref{ci-betti}.1)${}_n$
equality always holds for $n=0,1$, so the condition on $n$ cannot
be dropped from the last assertion.
 \end{Remark}

 \begin{proof}[Proof of Theorem \emph{\ref{ci-betti}}]
Let $\wh R\cong Q/\fa$ be a minimal Cohen presentation and set
$r=\rank_k\fa/\fn\fa$.  As $c+d=\edim R=\rank_k\pi^1_R$ and
$r=\rank_k\pi^2_R$, see \eqref{cohen}, from \eqref{lie-series} one
gets the first coefficient-wise inequality of formal power series
below:
 \begin{equation*}
\label{formal} \begin{aligned}
 \sum_{n=0}^\infty\rank_k\Ext nRkk\,t^n
&=\frac{(1+t)^{c+d}}{(1-t^2)^r}
  \frac{\prod_{i\ges1}(1+t^{2i+1})^{\rank_k\pi^{2i+1}}}
   {\prod_{i\ges1}(1-t^{2i+2})^{\rank_k\pi^{2i+2}}}\\
&\succcurlyeq\frac{(1+t)^{c+d}}{(1-t^2)^{r}}\\
&=\frac{(1+t)^{d}}{(1-t)^{c}}\frac{1}{(1-t^2)^{r-c}}\\
&\succcurlyeq\frac{(1+t)^d}{(1-t)^{c}}+(r-c)t^2\frac{(1+t)^d}{(1-t)^{c}}
 \end{aligned}
  \end{equation*}
The second inequality holds because $r-c$ is non-negative, see \eqref{codimension}.

Comparing coefficients, one obtains for every $n\ge0$ a numerical
inequality
 \[
\rank_k\Ext nRkk\ge \sum_{i=0}^{d}\binom
di\binom{c+n-i-1}{c-1}+(r-c)\sum_{i=0}^{d}\binom
di\binom{c+n-i-3}{c-1}\,,
 \]
which shows that the inequality in \ref{absolute_betti} holds for
every $n\ge0$.

If equality holds for some $n\ge2$, then the last formula yields
$r=c$, so by the Cohen-Macaulay Theorem the ideal $\fa$ is
generated by a regular sequence.

When $R$ is complete intersection one has $r=c$, and Tate
\cite[Thm. 6]{Ta} proves
 \[
\sum_{n=0}^\infty\rank_k\Ext nRkk
\,t^n=\frac{(1+t)^d}{(1-t)^{c}}\,.
 \]
This equality of power series means that \ref{absolute_betti}
holds for each $n\ge0$.
 \end{proof}

Regular local rings are precisely the complete intersection rings
of codimension $0$, and every complete intersection ring is
Gorenstein.  This hierarchy may also be observed by comparing 
the next result with Corollaries \ref{regular} and \ref{gorenstein}.

\begin{corollary}
\label{stable-ci-betti} For each $n\in\BZ$ there is an inequality
\begin{equation*}
\label{stable-betti} \tag*{\rm(\ref{stable-ci-betti}.1)${}_n$}
\rank_k\Extv nRkk\ge\sum_{i=0}^{d}\binom
di\left(\binom{c+n-i-1}{c-1} +\binom{c+d-n-i-2}{c-1}\right)\,.
 \end{equation*}

When $R$ is complete intersection equalities hold for each $n\in\BZ$.

If equality holds for a single $n\le d-3$ or $n\ge2$, then $R$ is
complete intersection.
 \end{corollary}

\setcounter{equation}{1}
\begin{Remark}
 \label{ci-remark2}
For $d\ge4$ there are no restrictions on $n$ in the last assertion,
so the value of $\rank_k\Extv nRkk$ for \emph{any} $n\in\BZ$
determines whether $R$ is complete intersection.  On the other hand, 
when $d$ satisfies $0\le d\le3$ Remark \ref{ci-remark1} and Corollary
\ref{gorenstein} show that for each $n$ in the non-empty interval 
$[d-2,1]$ the value of $\rank_k\Extv nRkk$ is the same for all 
Gorenstein rings with $\dim R=d$.
 \end{Remark}

\begin{proof}[Proof of Corollary \emph{\ref{stable-ci-betti}}] 
Corollary \ref{gorenstein} and Theorem \ref{ci-betti} show that inequalities
always hold in \ref{stable-betti}, and they become equalities when
$R$ is complete intersection.

If equality holds in \ref{stable-betti} for some $n\ge 2$ (respectively,
$n\le d-3$), then Corollary \ref{gorenstein} implies that equality
holds in \ref{absolute_betti} with $n\ge2$ (respectively, in
(\ref{ci-betti}.1)${}_{d-1-n}$ with $d-1-n\ge2$), so $R$ is complete
intersection by Theorem \ref{ci-betti}.
 \end{proof}

To put in context some results in this section we provide two examples.
The first one and Corollary \ref{stable-pd} shows that finiteness of
stable cohomology is not symmetric.

\begin{example}
 \label{arbitrary}
Let $(S,\fn,k)$ be an arbitrary local ring and set $R=S[x]/(x^2)$.
The $R$-module $M=R/(x)$ then has $\Extv nRMk\cong k$ for 
every $n\in\BZ$.

Indeed, the following sequence clearly is a complete resolution of
$M$ over $R$:
 \[
T=\quad\cdots\lra R\xra{\ x\ } R\xra{\ x\ }R\xra{\ x\ }R\lra\cdots
 \]
so Corollary \ref{complete-res} yields $\Extv nRMk\cong\CH
n{\Hom RTk}\cong k$ for each $n\in\BZ$.
 \end{example}

The next example should be compared to Corollary \ref{gorenstein}
and Theorem \ref{ci-betti}.

 \begin{example}
 \label{ambiguous}
Let $k$ be a field and $e$ a non-negative integer.  The ring
 \[
R=\frac{k[t_1,\cdots,t_e]}
 {\big(\{t_i^2-t^2_{i+1}\}_{1\le i\le e-1}
 \cup\{t_it_j\}_{1\le i<j\le e}\big)}
 \]
is artinian and Gorenstein for every $e\ge2$, the ring
 \[
S=\frac{k[t_1,\cdots,t_e]}
 {\big(\{t_1^2\}\cup\{t_1t_j\}_{3\le j\le e}
 \cup\{t_it_j\}_{2\le i\le j\le e}\big)}
 \]
is artinian, not Gorenstein for each $e\ge3$, and the following
equalities
 \[
\sum_{n=0}^\infty\rank_k\Ext n{R}kk\,t^n=\frac1{1-et+t^2}=
 \sum_{n=0}^\infty\rank_k\Ext n{S}kk\,t^n\,,
 \]
hold for all $e\ge2$, see \cite[Thm. 2]{LA} and \cite[Cor.\ p.\
38]{Fr}, respectively.
 \end{example}

\section{Structure of stable cohomology algebras}
\label{Structure of stable cohomology algebras}

This is the first of four sections devoted to the structure of the stable 
cohomology algebra of a local ring $(R,\fm,k)$.  We fix the following
notation.

\begin{chunk}
 \label{notation}
As in \eqref{sequence}, let $\iota$ be the canonical homomorphism
of graded $k$-algebras
 \[
\iota\col\CE\to\CS
 \quad\text{where}
\quad\CE=\Ext{}Rkk
 \quad\text{and}\quad
\CS=\Extv{}Rkk\,.
 \]

To describe the position of $\iota(\CE)$ in $\CS$ we use the
\emph{left torsion submodule} of $\CS$:
 \[
\CT=\{\sigma\in\CS\mid\CE^{\ges i}\cdot\sigma=0 \text{ for some
}i\ge0\}\,.
 \]
Note that one has $\CT=\gam\CS$, see  \eqref{gamma},
and that $\CT$ is an $\CE$-subbimodule of $\CS$.
 
Set $\CI=\Hom k{\CE}k$.  This is a graded $\CE$-bimodule with
the canonical actions:
 \[
 \begin{aligned}
(\varepsilon\cdot e)(\varepsilon')
&=(-1)^{|\varepsilon|(|e|+|\varepsilon'|)}e(\varepsilon'\cdot\varepsilon)
 \\
(e\cdot\varepsilon)(\varepsilon')
&=e(\varepsilon\cdot\varepsilon')
 \end{aligned}
\qquad\text{for all}\quad\varepsilon,\varepsilon'\in\CE
\quad\text{and}\quad e\in\CI\,.
 \]
The left action of $\CE$ on $\CI$ is of prime importance in later
developments.
 \end{chunk}

Regular rings are excluded from the theorem because for them
$\CS=0$, see \ref{regular}.

\begin{theorem}
 \label{structure}
Let $R$ be a singular local ring and set $\depth{}R=d$.
 \begin{enumerate}[\rm\quad(1)]
 \item
There is an exact sequence of left $\CE$-modules
 \begin{equation*}
 \label{general_sequence}
0\lra\CE\xra{\ \iota\ }\CS\xra{\ \eth\ }
\coprod_{i=d-1}^\infty(\shift^{-i}\CI)^{\mu^{i+1}}\lra0
  \end{equation*}
where $\mu^i=\rank_k\Ext iRkR$, and there are equalities
 \[
\CS=\iota(\CE)+\CE\cdot\CS^{\sles0}
 \quad\text{and}\quad
\iota(\CE)\cap\CT=0 \,.
 \]
 \item
If $\CS=\iota(\CE)\oplus\CT'$ for some graded left $\CE$-submodule
$\CT'\subseteq\CS$, then $\CT'=\CT$ and
 \[
\CT\cong\coprod_{i=d-1}^\infty(\shift^{-i}\CI)^{\mu^{i+1}}
 \quad\text{as graded left $\CE$-modules}\,.
 \]
 \item
If $\depth{}\CE\ge2$, then $\CS=\iota(\CE)\oplus\CT$ as graded
$\CE$-bimodules.
 \item
If $\wh R=Q/(f)$ for some singular local ring $(Q,\fn,k)$ and a
non-zero-divisor $f\in\fn^2$, then $\depth{}\CE\ge2$ and $\CT$ is
a two-sided ideal of $\CS$, such that
 \[
\CS=\iota(\CE)\oplus\CT\quad\text{and}\quad \CT\cdot\CT=0\,.
 \]
 \end{enumerate}
  \end{theorem}

For the proof we need a lemma.  It is well known and easy to show
that $\Ext{}Rkk$ and $\Tor{}Rkk$ are dual graded vector spaces.  A
more precise statement is contained in the next
result, which represents a version of Tate duality. It can be
derived, with a little work, from \cite[(2.1)]{Le}.  We provide a
direct argument.

\begin{lemma}
 \label{duality}
The graded left $\CE$-modules $\Tor{}Rkk$ and $\CI$ described in
\eqref{tor} and \eqref{notation}, respectively, have the
following properties.
 \begin{enumerate}[\quad\rm(1)]
  \item
There is a natural isomorphism $\delta\col \Tor{}Rkk\cong\CI$.
 \item
If $\vartheta$ is a right non-zero-divisor on $\CE$, then
$\CI=\vartheta^i\cdot\CI$ for each $i\ge0$.
\end{enumerate}
  \end{lemma}

\begin{proof}
(1)  Let $\varkappa\col G\to k$ be a minimal free resolution.

Setting $C=k$ in \eqref{tensor} one obtains the first map below:
\[
\Hom RGG\otimes_R(k\otimes_RG)\lra
k\otimes_RG\xra{{k\otimes_R\varkappa}}k\otimes_Rk=k\,.
\]
In homology the composition induces a morphism of graded $k$-vector
spaces
\[
\CE\otimes_k\Tor{}Rkk\lra k
\]
and hence a morphism of graded $k$-vector spaces
\begin{gather*}
\delta\col \Tor{}Rkk\lra\CI\quad\text{given by}\quad
\delta(1\otimes g)(\cls\alpha)=(-1)^{|g||\alpha|}\varkappa\alpha(g)
  \end{gather*}
A routine computation shows that $\delta$ is $\CE$-linear.

Any non-zero cycle in $k\otimes_RG_n$ has the form $1\otimes g$
for some $g\in G_n\smallsetminus\fm G_n$.  Choose an $R$-linear map
$\beta\col G_{n}\to k$ with $\beta(g)=1$. The Comparison Theorem
for resolutions yields a chain map $\alpha\col G\to G$ of degree
$-n$ with $\varkappa\alpha_n=\beta$; in particular,
\[
\delta({1\otimes g})(\cls\alpha)=\pm\beta(g)=\pm1\ne0\,.
\]
Thus, the $k$-linear map $\delta_n\col k\otimes_R G_n\to\CI^n$ is
injective.  Since the $k$-vector spaces $\Tor nRkk$ and $\Ext nRkk$
have the same finite rank, $\delta_n$ is even bijective.

(2) When $\vartheta$ is a right non-zero-divisor on $\CE$ the map
$\varepsilon\mapsto\varepsilon\cdot\vartheta^{i}$ is a $k$-linear
injection $\CE\to\CE$.  Its dual is a surjection $\CI\to\CI$ given
by $\upsilon\mapsto\pm\vartheta^i\cdot\upsilon$, so
$\CI=\vartheta^i\cdot\CI$.
 \end{proof}

 \begin{proof}[Proof of Theorem \emph{\ref{structure}}.]
(1)  Theorem \ref{martsinkovsky} and \eqref{sequence} yield an exact 
sequence
  \begin{equation}
 \label{exact}
0\lra\CE\xra{\ \iota\ }\CS\xra{\ \eth\ } \shift{\CB}\lra0
 \end{equation}
of graded $\CE$-bimodules with $\CB=\Extb{}Rkk$; the proof of the
theorem also shows
 \begin{equation}
 \label{basics}
{\CB}=\gam\CB
 \quad\text{and}\quad
\gam\CE=0\,.
 \end{equation}

Theorem \ref{bass-numbers} and Lemma \ref{duality}
yield isomorphisms of left $\CE$-modules
 \begin{equation}
  \label{cokernel}
 \begin{aligned}
 \shift^{-1}\Coker(\iota)&\cong\CB\cong\Ext{}RkR\otimes_k\Tor{}Rkk
\cong\Ext{}RkR\otimes_k\CI
\\
 &\cong\coprod_{i=d}^\infty(\shift^{-i}\CI)^{\mu^{i}}
\end{aligned}
 \end{equation}

{}From the definition of $\Gamma$ and \eqref{basics} we obtain
 \begin{equation}
 \label{intersection}
\iota(\CE)\cap\CT=\gam(\iota(\CE))\cong\gam\CE=0\,.
 \end{equation}

Next we set $\CN=\CE\cdot\CS^{\sles0}$ and prove
$\CS=\iota(\CE)+\CN$.  Note that for each $n<0$ the map $\eth^n$
yields $\CN^n=\CS^n\cong\CB^{n-1}$. Let $\vartheta\in\CE^2$ be a
right non-zero-divisor, see Lemma \ref{singular}. As $\shift\CB$
is a direct sum of shifts of $\CI$, Lemma \ref{duality}(2)
implies
 \begin{equation}
\label{divisible} \shift\CB=\vartheta^i\shift\CB
 \quad\text{for each}\quad i\ge0\,.
 \end{equation}

Let $\sigma$ be an arbitrary element of $\CS$.  By the remarks
above, one has $\eth(\sigma)=\vartheta^i\upsilon$ for some $\upsilon\in
(\shift\CB)^{\sles0}$ and $i\ge0$, and $\eth(\nu)=\upsilon$ for
some $\nu\in\CN$. Thus,
 \[
\eth(\sigma-\vartheta^i\nu)=\eth(\sigma)-\vartheta^i\upsilon=0\,,
 \]
hence $\sigma-\vartheta^i\nu$ is in $\Ker(\eth)=\iota(\CE)$, so one
gets $\sigma\in\vartheta^i\nu+\iota(\CE)\subseteq\CN+\iota(\CE)$.

(2)  By hypothesis, $\CS=\iota(\CE)\oplus\CT'$ for some left
$\CE$-submodule $\CT'\subseteq\CS$, so \eqref{exact} yields
$\CT'\cong\shift{\CB}$ . {}From \eqref{basics} we get
$\CT'=\gam{\CT'}\subseteq\gam{\CS}=\CT$, hence
$\iota(\CE)\cap\CT'=0$ by \eqref{intersection}. These relations
imply $\CT'=\CT$ and yield a direct sum of $\CE$-bimodules
 \begin{equation}
 \label{split}
\CS=\iota(\CE)\oplus\CT\,.
 \end{equation}
The expression for $\CT$ comes from the equality above and
the isomorphisms in \eqref{cokernel}.

(3)  If $\depth{}\CE\ge2$ holds, then one has
$\Ext1\CE{\shift\CB}\CE=0$, see Proposition \ref{torsion}, so
the sequence \eqref{exact} of graded left $\CE$-modules splits,
hence $\CS=\iota(\CE)\oplus\CT$, see (2).

(4) By Lemma \ref{residue-extension} we may assume $R=Q/(f)$, so
we obtain
 \[
\depth{}\CE=\depth{}{\Ext{}Qkk}+1\ge2
 \]
by referring to Proposition \ref{hypersection2} and Lemma
\ref{singular}. Part (3) now yields a direct sum decomposition
\eqref{split} of $\CE$-bimodules.  Next we prove $\CT\cdot\CT=0$.

Let $\vartheta\in\CE^2$ be the central non-zero-divisor defined by
the equality $R=Q/(f)$, see Proposition \ref{hypersection2}. Let
$\tau$ be an element of $\CT$.  Choosing $i\ge0$ so that
$\vartheta^i\cdot\tau=0$ holds, we get
$\tau\cdot\vartheta^i=\vartheta^i\cdot\tau=0$. On the other hand,
\eqref{exact} implies $\CT\cong\shift{\CB}$ as left $\CE$-modules,
so from \eqref{divisible} we get $\CT=\vartheta^{i}\cdot\CT$.
Thus, we obtain
 \[
\tau\cdot\CT=\tau\cdot(\vartheta^{i}\cdot\CT)=
(\tau\cdot\vartheta^{i})\cdot\CT=0\cdot\CT=0\,.
 \]
As $\tau\in\CT$ was chosen arbitrarily, we conclude
$\CT\cdot\CT=0$, as desired.

To finish the proof we apply the remark below.
 \end{proof}
 
 \begin{Remark}
 \label{two-sided-ideal}
If  $\CS=\iota(\CE)+\CT$ and $\CT\cdot\CT=0$, then 
$\CT$ is a two-sided ideal of $\CS$.

Indeed, $\CT$ is stable under multiplication on either side with elements 
of $\CE$, because it is an $\CE$-bimodule, and by its own elements,
as $\CT\cdot\CT=0$; thus, $\CS\cdot\CT\subseteq\CT\supseteq\CT\cdot\CS$.
 \end{Remark}

The following example provides applications for Theorem \ref{structure}(3).

 \begin{example}
 \label{local-tensor-product}
Let $k$ be a field, let $S$ and $T$ be localizations of finitely
generated $k$-algebras at $k$-rational maximal ideals. In this
case, $S\otimes_kT$ is a local ring.

Set $\CF=\Ext{}Skk$ and $\CG=\Ext{}Tkk$.  Standard K\"unneth
arguments give an isomorphism
$\Ext{}{S\otimes_kT}kk\cong\CF\otimes_k\CG$ of graded
$k$-algebras, so \cite[(3.1.iii)]{FHJLT} yields
 \[
\depth{}{\Ext{}{S\otimes_kT}kk}=\depth{}\CF+\depth{}\CG\,.
 \]
Thus, for singular $S$ and $T$ one has
$\depth{}{\Ext{}{S\otimes_kT}kk}\ge2$, see Lemma \ref{singular}.
  \end{example}

Applications of Theorem \ref{structure}(4) are discussed 
in the next section.

\section{Stable cohomology algebras of complete intersection rings}
\label{Stable cohomology algebras of complete intersection rings}

In this section we provide an explicit and nearly complete computation of 
the structure of the stable cohomology algebra of a complete 
intersection local ring.  This is achieved by combining results 
from this paper with already available information on the absolute 
cohomology algebra.   We fix notation for the entire section.

\begin{chunk}
 \label{presentation}
Let $(R,\fm,k)$ be a complete intersection local ring.  Fix a 
presentation $\wh R=Q/(\bsf)$ with $(Q,\fn,k)$ a regular 
local ring and $\bsf=\{f_1,\dots,f_c\}$ is a $Q$-regular sequence in $\fn^2$.
One then has $\edim R=e$, $\codim R=c$, and $\dim R=e-c$.

Fix a minimal generating set $\{t_1,\dots,t_e\}$ for $\fn$ and write
  \[
f_h\equiv\sum_{1\les i\les j\les e} a_{hij}\,t_it_j\pmod{\fn^3}
\quad\text{with}\quad a_{hij}\in Q
 \quad\text{for}\quad 1\le h\le c\,.
 \]
 
As in \eqref{notation}, we set $\CE=\Ext{}Rkk$ and $\CS=\Extv{}Rkk$.
 \end{chunk}

Sj\"odin \cite{Sj1} has completely described the algebra $\CE$; we
recall his result below.

\begin{chunk}
 \label{sjodin}
The algebra $\CE=\Ext{}Rkk$ is generated by elements
$\xi_1,\dots,\xi_e$ of degree $1$ and
$\vartheta_1,\dots,\vartheta_c$ of degree $2$, subject only to the
relations
  \begin{alignat*}{2}
\xi_i^2&=\sum_{h=1}^c \ov a_{hii}\vartheta_h
&&\quad\text{for}\quad 1\le i\le e\,;\\
\xi_i\xi_j+\xi_j\xi_i&=\sum_{h=1}^c\ov a_{hij}\vartheta_h
&&\quad\text{for}\quad 1\le i<j\le e\,;\\
\vartheta_h\xi_i&=\xi_i\vartheta_h &&\quad\text{for}\quad  1\le
h\le c \quad\text{and}\quad 1\le i\le e\,;\\
\vartheta_g\vartheta_h&=\vartheta_h\vartheta_g
&&\quad\text{for}\quad 1\le g\le h\le c\,,
 \end{alignat*}
where $\ov a_{hij}$ denotes the image of $a_{hij}$ in $k$; see
\cite[Thm.~5]{Sj1} or \cite[(10.2.2)]{res}.
 \end{chunk}

The next lemma determines the applicability of Theorem \ref{structure}(4).

\begin{lemma}
 \label{ci-depth}
For $R$ as in \emph{\eqref{presentation}} one has $\depth{}\CE=
\codim R$.
 \end{lemma}

 \begin{proof}
By Lemma \ref{residue-extension} we may assume $R=Q/(\bsf)$ for
a regular local $(Q,\fn,k)$ ring and a $Q$-regular set
$\bsf\subseteq\fn^2$, with $\operatorname{card}(\bsf)=\codim R$.
Proposition \ref{hypersection2} and Lemma \ref{singular} now
yield $\depth{}\CE=\depth{}{\Ext{}Qkk}+\codim R=\codim R$.
 \end{proof}

For complete intersection rings of codimension one, also known as 
\emph{hypersurface rings}, $\CS$ was computed by 
Buchweitz \cite[(10.2.3)]{Bu}.  We recover his result:

\begin{proposition}
 \label{hypersurface}
Let $R$ be a singular hypersurface ring of dimension $d$.

The $k$-algebra $\CS$ has generators $\xi_1,\dots,\xi_{e}$ of
degree $1$, a generator $\vartheta=\vartheta_1$ of degree $2$ and
a generator $\vartheta'$ of degree $-2$, subject to the relations
in \eqref{sjodin} and to
  \begin{gather*}
\vartheta\vartheta'=1= \vartheta'\vartheta\,.
 \end{gather*}
In particular, for each $n\in\BZ$ one has $\rank_k\CS^n=2^d$  
and $\Ext n{\CE}k{\CS}=0$.
  \end{proposition}

\begin{proof}
By Lemma \ref{residue-extension} we may assume $R\cong
Q/(f)$ with $(Q,\fn,k)$ regular local and
$f\in\fn^2\smallsetminus\{0\}$. By Corollary \ref{localization} the
algebra $\CS$ is the localization of $\CE$ at $\{\vartheta^s\var
s\ge0\}$, so the presentation of $\CE$ in \eqref{sjodin} yields
the presentation above.

The $k$-rank of $\CS^n$ can be obtained from this presentation, or 
directly from Corollary \ref{stable-ci-betti}.  As $\vartheta$ is a central 
element in $\CE$ and is invertible in $\CS$, \eqref{rees}(1) yields
 \[
\depth\CE\CS=\depth\CE{(\CS/\CS\vartheta)}+1=\depth\CE0+1=\infty\,.
\qedhere
 \]
 \end{proof}
 
In higher codimension the structure of stable cohomology 
is completely different. 

\begin{proposition}
 \label{ci2}
Let $R$ be a complete intersection ring with $\codim R=c\ge2$.

The graded $k$-algebra $\CS$ has the form $\CS\cong\CE\oplus\CT$, 
where $\CE$ is the graded algebra from \eqref{sjodin}, 
$\CT=\bigcup_{i\ges 1}\{\tau\in\CS\var\CS^{\ges i}\tau=0\}$ 
is an ideal with $\CT\cdot\CT=0$, $\CT\cong\Hom k\CE k$
as graded left $\CE$-module, and $\tau\cdot\vartheta_h=
\vartheta_h\cdot\tau$ for all $\tau\in\CT$.
 \end{proposition}
 
\begin{proof}
Theorem \ref{structure}(4) gives everything but the last assertion, which
amounts to saying that each element $\iota(\vartheta_h)$ is central in $\CS$.  
To see this note that $\vartheta_h\in\CE^2$ can be defined by the canonical 
surjection $Q/(\bsf\smallsetminus\{f_h\})\to R$ and apply Theorem
 \ref{central}.
  \end{proof}

Still missing from Proposition \ref{ci2} for a full description of the structure 
of $\CS$ are data on the products $\tau\cdot\xi_j$.  Specifically, we ask:

  \begin{question}
Is $\CT$ isomorphic to $\shift^{1-d}\Hom k\CE k$ as graded $\CE$-bimodules?
  \end{question}
  
Here is a case where a positive answer is available from other sources.

\begin{example}
 \label{groups}
Let $k$ be a field of characteristic $p>0$, let $u_1,\dots,u_e$ be positive
integers, and set $q_i=p^{u_i}$.  For $R=k[x_1,\dots,x_e]/(x_1^{q_e},\dots,
x_1^{q_1})$ the algebra $\CS$ is graded-commutative, so one has 
$\CT\cong \Hom k\CE k$ as graded $\CE$-bimodules.

Indeed, one has $R\cong kG$, where $G$ is the abelian group
$\prod_{i=1}^e\BZ/(p^{u_i})$, and compatible isomorphisms of graded 
$k$-algebras $\CE\cong\CH*{G,k}$ and $\CS\cong\wh{\operatorname H}{}^*{(G,k)}$,
with the ordinary cohomology and the Tate cohomology of $G$; see \cite[(V.4.6)]{Br}, 
\cite[(VI.6.2)]{Br} and \cite[(6.11)]{BKS}.  The algebra 
$\wh{\operatorname H}{}^*{(G,k)}$ is graded-commutative; see \cite[(XII.5.3)]{CE}.
 \end{example}

\section{Stable cohomology algebras of Gorenstein rings}
\label{Stable cohomology over Gorenstein rings}

In this section $(R,\fm,k)$ is a Gorenstein local ring of
dimension $d$.  

The leading theme here would come as no surprise: stable cohomology 
is simpler over Gorenstein rings.  Numerically, this has already been made
precise by Corollary \ref{gorenstein}.  We back it up by showing that
the structure of the stable cohomology algebra of a Gorenstein ring is
much more rigid than might have been expected \emph{a priori}.  Although 
not as explicit as the in the special case of complete intersections, treated
in Section \ref{Stable cohomology algebras of complete intersection
rings}, the results here are significantly  sharper than 
those in Section \ref{Structure of stable cohomology algebras}. 

The following notation and terminology is in force for the entire section.

\begin{chunk}
 \label{split2}
Set $\CE=\Ext{}Rkk$ and $\CS=\Extv{}Rkk$, let $\iota\col\CE\to\CS$
denote the canonical homomorphism of graded $k$-algebras, $\CT$ the
torsion subbimodule $\gam\CS\subseteq\CS$ and $\CI$
the graded left $\CE$-module $\Hom k{\CE}k$ with the natural action, 
see \eqref{notation}

We say that the algebra $\CS$ \emph{splits} if $\CT$ is a two-sided
ideal of $\CS$, such that 
 \[
\CT\cdot\CT=0\,,\quad
\CS=\iota(\CE)\oplus\CT\,,
\quad\text{and}\quad
\CT\cong\shift^{1-d}\CI\quad\text{as graded left $\CE$-modules}\,.
 \]
 \end{chunk}

The structure of $\CS$ for Gorenstein rings $R$ with $\codim R\le
1$ is completely known, see Proposition \ref{hypersurface}, so this
case is excluded from the next theorem.

\begin{theorem}
 \label{structure2}
Let $R$ be a Gorenstein local ring with $\dim R=d$ and $\codim
R\ge2$.
 \begin{enumerate}[\rm\quad(1)]
 \item
There is an exact sequence of left $\CE$-modules
 \begin{equation*}
   \label{gor_sequence}
0\lra\CE\xra{\ \iota\ }\CS\xra{\ \eth\ } \shift^{1-d}\CI\lra0
  \end{equation*}
and there are equalities $\CS=\iota(\CE)+\CE\cdot\CS^{\sles0}$ and
$\iota(\CE)\cap\CT=0$.
 \item
If $\CS=\iota(\CE)\oplus\CT'$ for some graded left $\CE$-submodule
$\CT'\subseteq\CS$, then  $\CT'=\CT$ and the algebra $\CS$ splits.
 \item
If $\depth{}\CE\ge2$, then the algebra $\CS$ splits.
 \item
If $\zeta\cdot\CE=\CE\cdot\zeta$ for some left non-zero-divisor
$\zeta\in\CE^{\ges1}$, then $\CS$ splits and
 \[
\CT=\CE\cdot\CS^{\sles0}\, .
 \]
 \end{enumerate}
  \end{theorem}

The proof is presented after a few remarks.

\begin{Remark}
 \label{gor-remark1}
The theorem above should be compared to Theorem \ref{structure}, 
applied to a Gorenstein ring $R$.  Part (1) is a simple specialization.  
On the other hand,  in Theorem \ref{structure2} the conclusions of parts 
(2) through (4) are stronger, while the hypothesis of (4) is weaker, 
see Proposition \ref{hypersection2}.
 \end{Remark}

\begin{Remark}
 \label{gor-remark2}
Theorem \ref{structure2} offers striking parallels to results of
Benson and Carlson \cite{BC} relating the structure of the Tate
cohomology algebra $\wh{\rm H}{}^*(G,k)$ of a finite group $G$ to
that of the absolute cohomology algebra $\CH{*}{G,k}$.

Such similarities are unexpected, as the corresponding 
algebras have completely different properties.
Indeed, $\CH{*}{G,k}$ and $\wh{\rm H}^*(G,k)$ are graded
commutative algebras, and the first one is also finitely generated
over $\HH{0}{G,k}=k$.  In stark contrast,
the algebra $\Ext{}Rkk$ may not be finitely generated; it is
noetherian if and only if $R$ is complete intersection; it is
commutative if and only if $\wh R\cong Q/\fa$, with $(Q,\fn,k)$
regular and $\fa$ generated by a $Q$-regular sequence contained in
$\fn^3$.
 \end{Remark}

\begin{Remark}
 \label{gor-remark3}
Parts (3) and (4) of Theorem \ref{structure2} do  not cover all cases when
the algebra $\CS$ splits; see  Example~\ref{felix2} and Question~\ref{split?}.  
 \end{Remark}

 \begin{proof}[Proof of Theorem \emph{\ref{structure2}}]
(2)  If $\CS=\iota(\CE)\oplus\CT'$ as left $\CE$-modules, then Theorem
\ref{structure}(2) yields an equality $\CT'=\CT$ and a decomposition
$\CS=\iota(\CE)\oplus\CT$ as $\CE$-bimodules.

To prove $\CT\cdot\CT=0$ we fix a right non-zero divisor
$\vartheta\in \CE^2$, see Lemma \ref{singular}.  Let $\tau$ be
an element of $\CT$. As $\CT$ is an $\CE$-subbimodule of $\CS$, 
and for $j\ge d$ one has $\CT^{j}\cong\CI^{j+1-d}=0$, we get
$\tau\cdot\vartheta^i\in\CT^{\lceil\tau\rceil+2i}=0$ for all $i\gg0$. Fix 
an integer $i$ with this property.  Lemma
\ref{duality}(2) yields $\CI=\vartheta^i\cdot\CI$, so (1) implies
$\CT=\vartheta^i\cdot\CT$, and hence
 \[
\tau\cdot\CT=\tau\cdot(\vartheta^{i}\cdot\CT)=
(\tau\cdot\vartheta^{i})\cdot\CT=0\cdot\CT=0\,.
 \]
As $\tau\in\CT$ was arbitrary, this implies $\CT\cdot\CT=0$,
so $\CT$ is an ideal by Remark \ref{two-sided-ideal}.

(3) This follows from (2), in view of Theorem \ref{structure}(3).

(4) Set $\CN=\CE\cdot\CS^{\sles0}$.  As $\iota(\CE)+\CN=\CS$ by
Theorem \ref{structure}(1), it suffices to prove
$\iota(\CE)\cap\CN=0$, see (2).  Assuming the contrary means
that one can write
 \[
0\ne\sum_{h=1}^u\varepsilon_h\cdot\nu_h\in\iota(\CE)
 \]
with $\varepsilon_h\in\CE$, $\nu_h\in\CS^{\sles0}$, and 
$\varepsilon_h\cdot\nu_h\ne0$.
Set $s=\lceil\zeta\rceil$ and choose $i\ge0$ so that
 \[
is+\lceil\nu_h\rceil>d \quad\text{holds for}\quad h=1,\dots,u\,.
 \]
As $\zeta$ is a non-zero-divisor on $\CE$ and $\zeta\cdot\CE=
\CE\cdot\zeta$, there exist $\varepsilon'_h\in\CE$ such that
 \[
0\ne\sum_{h=1}^u\zeta^i\cdot\varepsilon_h\cdot\nu_h=
\sum_{h=1}^u\varepsilon'_h\cdot\zeta^i\cdot\nu_h\,.
 \]

Choose $l\in[1,u]$ with $\varepsilon'_l\cdot\zeta^i\cdot\nu_l\ne0$, then set 
$\nu=\nu_l$ and $n=\lceil\nu\rceil$; one then has
 \[
0\ne\zeta^i\cdot\nu\in\CS^{is+n}\subseteq\CS^{\sges d} =\iota(\CE)^{\sges
d}\subseteq\iota(\CE)^{\ges1}\,.
 \]
The map $\eth$ in (1) induces a bijection 
$\CS^{\sles0}\to (\shift^{1-d}\CI)^{\sles0}$.  Set $r=\rank_k\CN^n$.
As $\CI=\zeta^{2r}\CI$ by Lemma \ref{duality}(2), we may choose 
$\nu'\in\CN^{n-2rs}$ with $\zeta^{2r}\cdot\nu'=\nu$.  

Let $\cD$ be the universal enveloping algebra of the graded Lie
subalgebra $\pi^{\operatorname{even}}_R$ of $\pi_R$.  By
\eqref{PBW-nzd}, it has no zero-divisors 
(different from $0$) and $\CE$ is a free graded $\cD$-module.  
Thus, $\cD^{2j}$ and $\cD^{2j}\cdot\zeta^i\cdot\nu$ are isomorphic 
$k$-spaces for each $j\in\BZ$. From
 \begin{align*}
\cD^{2rs} &\cong\cD^{2rs}\cdot\zeta^i\cdot\nu_1
=\cD^{2rs}\cdot\zeta^{i+2r}\cdot\nu'\\
&\subseteq\CE^{2rs}\cdot\zeta^{i+2r}\cdot\nu'
=\zeta^{i+2r}\cdot\CE^{2rs}\cdot\nu'\\
&\subseteq\zeta^{i+2r}\cdot\CN^{n}
 \end{align*}
we get $\rank_k\cD^{2rs}\le\rank_k\CN^n=r$.  On the other hand,
we have $r_2\ge\codim R\ge2$, see \eqref{codimension},
so\eqref{lie-series} yields coefficient-wise inequalities of formal
power series
 \begin{align*}
\sum_{j=0}^\infty(\rank_k\cD^{2j})t^{2j} &=\frac1{(1-t^2)^2}\cdot
\frac1{(1-t^2)^{r_{2}-2}\prod_{j\ges2}(1-t^{2j})^{r_{2j}}}
\\
&\succcurlyeq\frac1{(1-t^2)^2}=\sum_{j=0}^{\infty}(j+1)t^{2j}
 \end{align*}
In particular, we get $\rank_k\cD^{2rs}>r$, and
hence the desired contradiction.
 \end{proof}

We proceed to test the hypotheses of the preceding theorem.
In the rest of this section the discussion revolves around the
following question:

 \begin{question}
  \label{split?}
Does the cohomology algebra $\CS=\Extv {}Rkk$ split for every
Gorenstein local ring $(R,\fm,k)$ with $\codim R\ge 2$?
 \end{question}

We start with a case when Theorem \ref{structure2}(4) applies.

 \begin{example}
 \label{conormal}
If $\wh R\cong Q/\fa$ is a minimal Cohen presentation
and the $\wh R$-module $\fa/\fa^2$ has a non-zero free direct summand,
then Iyengar \cite[Main Theorem]{Iy} provides a non-zero central
element in $\pi^2_R$, hence a central non-zero-divisor in $\CE$,
cf.\ \eqref{PBW-nzd}.

It might be noted that $\fa/\fa^2$ has a non-zero free direct summand
whenever $\wh R$ is isomorphic to $Q/(f)$ for some local ring $(Q,\fn,k)$
and non-zero-divisor $f\in\fn^2$.  This case is covered already by Theorem
\ref{structure}(4), but it is not known whether all free direct summands
of $\fa/\fa^2$ arise in this way.
  \end{example}

The first application of Theorem \ref{structure2}(3) mirrors Example
\ref{local-tensor-product}.

 \begin{example}
 \label{tensor-product:G}
Let $k$ be a field, and let $S$ and $T$ be localizations of
finitely generated $k$-algebras at $k$-rational maximal ideals. If
$S$ and $T$ are singular Gorenstein rings, then so is the local
ring $S\otimes_kT$, and Example \ref{local-tensor-product} yields
 \[
\depth{}{\Ext{}{S\otimes_kT}kk}\ge2\,.
 \]
  \end{example}

To exhibit further classes of Gorenstein rings satisfying the hypothesis
of Theorem \ref{structure2}(3) we use Koszul duality, see \cite[\S 5]{Sm}
for a concise introduction.

 \begin{chunk}
 \label{artin_koszul}
Let $k$ be a field, and let $A=\bigoplus_{i\ges0}A_i$ be a commutative
internally graded $k$-algebra, see \eqref{gradings}, with $A_0=k$,
$A=k[A_1]$, and $\rank_kA_1<\infty$.

Recall that $A$ is said to be \emph{Koszul} if $\Ext nAkk^i=0$ for $n\ne
i$. Its \emph{Koszul dual} is the internally graded
$k$-algebra $B$ with $B_i=\Ext iAkk^i$ and multiplication induced by
the composition product of $\Ext{}Akk$; the algebra $B$ is Koszul
as well, and its own Koszul dual is canonically isomorphic to $A$.

If $A$ is a Gorenstein Koszul algebra with $\rank_kA<\infty$, then
 \[
\Ext nBkB\cong\begin{cases}
k&\text{for } n=\sup\{n\in\BN\mid A_n\ne0\}\,;\\
0&\text{otherwise}\,.
 \end{cases}
 \]

Indeed, this follows by Koszul duality from \cite[(5.10),
(4.3.1)]{Sm}.
  \end{chunk}

\begin{proposition}
 \label{gor_koszul}
Let $R$ be the localization of a commutative Gorenstein graded
$k$-algebra $A=\bigoplus_{i\ges0}A_i$ at the maximal ideal
$A_+=\bigoplus_{i\ges1}A_i$.  If $A$ is Koszul, then
 \[
\depth{}\CE\ge2
 \]
holds, unless $A\cong k[x_1,\dots,x_e]/(f)$ for some quadratic form $f$
(possibly equal to $0$).
  \end{proposition}

 \begin{proof}
Assuming $\depth{}\CE\le1$, we induce on $\dim A$ to prove that $A$
has the desired special form above.  When $\dim A=0$ one has $\rank_kA<\infty$,
so \eqref{artin_koszul} yields $A_i=0$ for $i\ge2$. As $A$ is Gorenstein,
this is only possible if $A\cong k$ or $A\cong k[x]/(x^2)$.

Suppose $\dim A=d\ge1$ and the assertion holds for algebras
of dimension $d-1$.  Let $k\to k'$ be a field extension. The
$k'$-algebra $A'=k'\otimes_kA$ is clearly Koszul. Let $R'$ be its
localization at $A'_+$. The canonical map $R\to R'$ is flat and
$R'/\fm R'\cong k'$, so Lemma \ref{residue-extension} yields
$\depth{}{\Ext{}{R'}{k'}{k'}}=\depth{}\CE$.  Thus, we may
assume $k$ is infinite, and so find a
non-zero-divisor $g\in A_1$.  Now $\ov A=A/(g)$ is Koszul with $\dim\ov
A=d-1$. Since $g/1\in R$ is $R$-regular and not in $\fm^2$, for $\ov
R=R/(g/1)$ Proposition \ref{hypersection1} yields $\depth{}{\Ext{}{\ov
R}kk}=\depth{}\CE$. As $\ov R$ is the localization of $\ov A$ at
$\bigoplus_{i\ges1}\ov A_i$, the induction hypothesis yields $\ov A\cong
k[\ov x_1,\dots,\ov x_{e-1}]/(\ov f)$ for some quadratic form $\ov f$.
It follows that $A$ has the desired property.
 \end{proof}

Let $\mult R$ denote the multiplicity of $R$.
A Gorenstein ring $R$ has multiplicity $1$ if and only if it is
regular, and multiplicity $2$ if and only if it is a quadratic
hypersurface; else, it satisfies $\mult R\ge\codim R+2$.

The case of minimal multiplicity is covered by the next result.

\begin{proposition}
 \label{minimal:G}
If $R$ is Gorenstein and $\mult R=\codim R+2$, then
 \[
\depth{}\CE=
 \begin{cases}
1&\text{when $R$ is a hypersurface ring}\,;
\\
 2&\text{otherwise}\,.
\end{cases}
 \]
  \end{proposition}

 \begin{proof}
One has $\codim R=1$ if and only if $R$ is a hypersurface, and when
this is the case Lemma \ref{ci-depth} applies.  For the rest of the
proof assume $\codim R\ge2$.  In view of Lemma \ref{weak} we
may further assume $\dim{}R=0$. Our hypothesis then implies
isomorphisms $\fm^3=0$ and $\fm^2=\ann_R\fm\cong k$, and 
non-degeneracy of the map
 \[
\mu\col(\fm/\fm^2)\otimes_k(\fm/\fm^2)\lra \fm^2\cong k
 \]
induced by the product of $R$.  Thus, $A=k\oplus(\fm/\fm^2)\oplus\fm^2$ 
is a Gorenstein graded $k$-algebra, with
$\rank_kA_1=\codim R$.  By \cite[(3.1)]{LA} both $\CE$ and $\Ext{}Akk$
are isomorphic to the tensor algebra of $(A_1)^\vee$ over $k$, modulo
the quadratic form $\Hom k\mu k(\id k)\in(A_1)^\vee\otimes_k(A_1)^\vee$.
Thus, we may replace $R$ by $A$. The description of $\Ext{}Akk$ shows
that as an algebra over $\Ext0Akk=k$ it is generated by $\Ext1Akk$.
This graded vector space has $\Ext1Akk^j=0$ for $j\ne1$, so $\Ext
iAkk^j=0$ for $j\ne i$, so $A$ is Koszul, hence $\depth{}\CE=2$ by
\eqref{artin_koszul}.
 \end{proof}

Recall that if $R$ is Gorenstein with $\codim R\le 2$, then $R$ is
complete intersection, so Lemma \ref{ci-depth} gives $\depth{}\CE=\codim
R$.  In the next codimension we prove:

\begin{proposition}
 \label{codim3}
If $R$ is Gorenstein and $\codim R=3$, then
 \[
\depth{}\CE=\begin{cases}
3&\text{when $R$ is complete intersection}\,;\\
2&\text{otherwise}\,.
\end{cases}
 \]
 \end{proposition}

 \begin{proof}
By Proposition \ref{ci-depth} we may assume $R$ is not complete
intersection.  By Lemma \ref{subalgebra} it suffices to prove
$\depth{}\cD=2$, where $\cD$ is the universal enveloping subalgebra of
$\pi^{\ges2}_R$.  It follows from \cite[(3.3)]{small} that $\pi^{\ges2}_R$
is isomorphic to $\pi_A$, where $A=\bigoplus_{i=0}^3\HH iK$ and $K$ is
the Koszul complex on a minimal set of generators of $\fm$.  There exist
bases $\{e_1,\dots,e_r\}$ of $A_1$; $\{f_1,\dots,f_r\}$ of $A_2$; $\{g\}$
of $A_3$, such that
 \[
e_if_i=g=f_ie_i
 \quad\text{for}\quad 1\le i\le r
  \]
and all other products of basis elements vanish, see \cite[(8.4)]{small}.
Thus, $A$ is a Gorenstein local ring of minimal multiplicity, so
Proposition \ref{minimal:G} yields $\depth{}\cD=2$.
 \end{proof}

In the balance of this section we deal with artinian rings.  

 \begin{proposition}
  \label{artinian}
For an artinian Gorenstein local ring $(R,\fm,k)$ with $\fm\ne0$
the following conditions are equivalent.
 \begin{enumerate}[\rm(i)]
 \item
$\CS$ is split.
 \item
$\CT=\CS^{\sles0}$.
 \item
$\CS^{\sles0}$ is a left $\CE$-submodule of $\CS$.
 \item
$\CE$ is generated as an algebra over $\CE^0=k$ by a set
$E\subseteq\CE^{\ges1}$, such that $\varepsilon\cdot\CS^j=0$
holds for all pairs $(\varepsilon,j)\in E\times\BZ$ satisfying
$-\lceil\varepsilon\rceil\le j<0$.
 \end{enumerate}
\end{proposition}

 \begin{proof}
{}From the exact sequence in Theorem \ref{structure2}(1), or from
Theorem \ref{perpendicular}, one gets $\iota(\CE)=\CS^{\ges0}$, hence
$\iota(\CE)$ has a unique complement in $\CS$ as a graded $k$-vector
space, namely, $\CS^{\sles0}$.  Thus, when $\CS$ is split one has
$\CT=\CS^{\sles0}$, and hence (i) implies (ii).

It is clear that (ii) implies (iii).  If $\CS^{\sles0}$ is a left
$\CE$-submodule of $\CS$, then it is necessarily a direct
complement of $\iota(\CE)$, so (iii) implies (i) by Theorem
\ref{structure2}(2).

It is clear that (iii) implies the validity of (iv) for \emph{every}
set of generators $E\subseteq\CE^{\ges1}$, in particular, for
$E=\CE^{\ges1}$.  Conversely, assume that (iv) holds for some $E$,
pick an arbitrary $\sigma\in\CS^{\sles0}$, and 
$\varepsilon_1,\dots,\varepsilon_s$ be elements of $E$.
If one has  $\lceil\sigma\rceil<-\sum_{i=1}^s\lceil\varepsilon_i\rceil$, 
then $\varepsilon_1\cdots\varepsilon_s\cdot\sigma\in\CS^{\sles0}$ 
holds for degree reasons.  Else, there is an integer $r$ with
 \[
1\le r\le s
 \quad\text{and}\quad
-\sum_{i=r}^s\lceil\varepsilon_i\rceil
\le\lceil\sigma\rceil<-\sum_{i=r+1}^s\lceil\varepsilon_i\rceil\,.
 \]
Thus, $\sigma'=\varepsilon_{r+1}\cdots\varepsilon_s\cdot\sigma$ satisfies
 $-\lceil\varepsilon_{r}\rceil\le\lceil\sigma'\rceil<0$, hence
 \[
\varepsilon_1\cdots\varepsilon_s\cdot\sigma=
(\varepsilon_1\cdots\varepsilon_{r-1})\cdot(\varepsilon_{r}\cdot\sigma')=0\,.
 \]
It follows that $\CS^{\sles0}$ is a left $\CE$-submodule of $\CS$,
that is, (iii) holds.
 \end{proof}

When $\dim R=0$ it is easy to write down a complete
resolution of $k$, see \ref{complete-res-def}.

 \begin{chunk}
  \label{artinian-resolution}
Let $(R,\fm,k)$ be an artinian Gorenstein ring with $\fm\ne0$.

If $\epsilon\col F\to k$ is a minimal free resolution and $(-)^*=\Hom
R-R$, then
 \[
T=\quad\cdots\lra
F_1\xra{\,\dd_1\,}F_0\xra{\,\epsilon_0^*\circ\,\epsilon_0\,}F_0^*
\xra{\dd_1^*}F_1^*\lra\cdots
 \]
where $T_0=F_0$, is a complete resolution of $k$ satisfying
$\dd(T)\subseteq\fm T$.

Indeed, one has $\dd(T)\subseteq\fm T$ by construction.  The exactness of
$T$ follows easily from the self-injectivity of $R$, which then implies
the exactness of $\Hom RTR$.
 \end{chunk}

We finish with an example obtained through the dictionary 
between local algebra and rational homotopy theory; see \cite{AH}.
For a topological space $Y$ let $\CH n{Y;\BQ}$ and $\HH n{Y;\BQ}$
denote its singular (co)ho\-mo\-logy, $\Omega Y$ its loop
space, and $\HH*{\Omega Y;\BQ}$ its Pontryagin algebra with
multiplication induced by composition of loops.  

The construction below was communicated to us by Yves F\'elix, in answer 
to our question whether there exists a formal CW complex $Y$ with
$\depth{}{\HH*{\Omega Y;\BQ}}=1$.

 \begin{chunk}
  \label{felix1}
Let $S^2$ denote the standard $2$-sphere, $\#$ a connected sum of
smooth manifolds.  The following manifolds are formal topological spaces:
 \[
X=S^2\times S^2\times S^2=X' \quad\text{and}\quad Y=X\#X'
 \]

Using a suitable CW decomposition of $X\#X'$, F\'elix exhibits 
$\HH*{\Omega Y;\BQ}$ as a free product of three commutative 
polynomial algebras:
 \begin{equation}
  \label{free-product}
\HH*{\Omega Y;\BQ}\cong
\BQ[\xi_1,\xi_2,\xi_3]*\BQ[\xi'_1,\xi'_2,\xi'_3]*\BQ[\vartheta]\,.
 \end{equation}
where $\HH*{\Omega X;\BQ}=\BQ[\xi_1,\xi_2,\xi_3]$, and
$\HH*{\Omega X';\BQ}=\BQ[\xi'_1,\xi'_2,\xi'_3]$ with
$\lfloor\xi_i\rfloor=1=\lfloor\xi'_j\rfloor$, and $\vartheta\in\HH4{\Omega Y;\BQ}$ 
arises from the identification in $Y=X\#X'$ of the orientations of $X$ and $X'$.
None of these algebras equals $\BQ$, so \cite[(36.e.2)]{FHT} yields
 \begin{equation}
  \label{free-depth}
\depth{}{\HH*{\Omega Y;\BQ}}=1\,.
 \end{equation}
   \end{chunk}

Next we translate Felix's example into commutative algebra.

 \begin{example}
 \label{felix2}
Let $k$ be a field of characteristic $0$.  The ring
 \[
R=\frac{k[t_1,t_2,t_3,t'_1,t'_2,t'_3]}
{\big(\{t_i^2\,,\,t_it'_j\,,\,t'^2_j\}_{1\le i,j\le 3}\,,\,
t_1t_2t_3-t'_1t'_2t'_3\big)}
 \]
is Gorenstein with $\codim R=6$, $\fm^3\ne0=\fm^4$, and
$\rank_kR=14$.

Its cohomology algebra $\CE$ satisfies $\depth{}\CE=1$ and
$\CE\zeta\ne\zeta\CE$ for every non-zero $\zeta\in\CE^{\ges1}$,
and its stable cohomology algebra $\CS$ is split.

Indeed, the properties of $R$ are clear. By Lemma
\ref{residue-extension} we may assume $k=\BQ$.

Endow $R$ with an internal grading by assigning (lower) degree $-2$ to the
variables $t_i$ and $t'_j$ one gets an isomorphism
$R\cong\bigoplus_{n=0}^6\CH n{Y;\BQ}$ of internally graded rings,
where $Y$ is the manifold from \eqref{felix1}. Since $Y$ is
formal, there are isomorphisms
 \[
\bigoplus_{l+m=-n}\Ext l{R}{\BQ}{\BQ}^{m} \cong\HH n{\Omega Y;\BQ}
 \]
of $\BQ$-vector spaces that combine into an isomorphism of graded
$\BQ$-algebras
 \[
\Ext{}{R}{\BQ}{\BQ}\cong\HH*{\Omega Y;\BQ}\,.
 \]
Thus, one obtains $\depth{}\CE=1$ from \eqref{free-depth}.

To analyze $\CS$ we note that a minimal graded free resolution $F$
of $\BQ$ starts as
 \[
\cdots\lra R(4)^{\binom{3+3}2}\oplus R(4)^{3+9+3}\oplus R(6)\lra
R(2)^{3+3}\lra R\lra0
 \]
where the standard basis of $R(2)^{3+3}$ maps to the generators
$t_1,t_2,t_3$ and $t'_1,t'_2,t'_3$ of $\fm$, that of
$R(4)^{\binom{3+3}2}$ to the Koszul relations between these
generators, that of $R(4)^{3+9+3}$ to the syzygies defined by the
quadratic relations of $R$, and that of $R(6)$ to the syzygy
defined by the cubic relation.  The complete resolution $T$
constructed in \eqref{artinian-resolution} from the minimal
resolution $F$ now yields isomorphisms
 \begin{alignat*}{4}
\CE^2=&&\CS^2&\cong \BQ(4)^{30}\oplus \BQ(6)
\qquad \CE^0=&\CS^0&\cong \BQ
\qquad &\CS^{-1}&\cong \BQ(-6)
 \\
\CE^1=&&\CS^1&\cong \BQ(2)^{6}
&& &\CS^{-2}&\cong \BQ(-8)^6
 \end{alignat*}
The isomorphism \eqref{free-product} shows that the algebra $\CE$ is
generated over $\CE^0=\BQ$ by elements from $\CE^1$ and $\CE^2$.
As the action of $\CE$ on $\CS$ is compatible with internal
gradings, see Proposition \ref{stable-graded}, degree
considerations yield the following equalities:
 \[
\CE^1\cdot\CS^{-1}=0\,,\quad\CE^2\cdot\CS^{-1}=0\,,
\quad\CE^2\cdot\CS^{-2}=0\,.
 \]
Proposition \ref{artinian} now shows that the stable cohomology
algebra $\CS$ is split.
 \end{example}

\section{Stable cohomology algebras of Golod rings}
\label{Stable cohomology over Golod rings}

In this section $(R,\fm,k)$ denotes a local ring.  The \emph{codepth}
of $R$ is the integer $\codepth R=\edim R-\depth{}R$.  Once again,
we consider the graded $k$-algebras
 \[
\CE=\Ext{}Rkk
 \quad\text{and}\quad
\CS=\Extv{}Rkk
 \]
and let $\iota\col\CE\to\CS$ denote the canonical homomorphism of 
graded algebras.

For all rings $R$ with $\codepth R\ge2$ analyzed so 
far in this paper, $\iota(\CE)$ has a direct complement in $\CS$ 
as a left $\CE$-submodule.  Here our goal is to produce a ring 
$R$ for which this fails.  We search for it among Golod rings, 
as their homological properties are in many respects
antithetical to those of Gorenstein rings. 

Golod rings are usually defined in terms of the series 
$\sum_{n\ges 0}\rank_k\CE^nt^n$, see \cite[\S 5]{res} for 
details and examples. Here it is useful to take as definition their 
characterization in terms of the structure of the graded $k$-algebra $\CE$.

\begin{chunk}
 \label{golod-def}
The ring $R$ is said to be \emph{Golod} if the universal
enveloping algebra $\cD$ of the Lie algebra $\pi^{\ges2}_R$ is a
free associative $k$-algebra, \cite[Cor.\ p.\ 59]{gol}.
 \end{chunk}

We will also use a defining homological property of free $k$-algebras.

\begin{chunk}
  \label{tensor-algebra}
Let $\CA$ be a free associative $k$-algebra on a set $\Xi$ of
elements of positive upper degree. Let $\{b_\xi\var\lceil
b_\xi\rceil=\lceil\xi\rceil+1\}_{\xi\in\Xi}$ be a linearly
independent set over $\CA$.  There is then an exact sequence of
graded left $\CA$-modules:
 \begin{equation}
 0\lra\coprod_{\xi\in\Xi}\CA b_\xi\xra{\ \dd\ }\CA\lra k\lra0
\quad\text{with}\quad
 \dd(b_\xi)=\xi\,.
  \end{equation}
 \end{chunk}

It is clear that $R$ is Golod when $\codepth R=0$ (because then
$R$ is regular, so $\cD=k$), or when $\codepth R=1$ (because then
$R$  is a singular hypersurface ring, so $\cD=k[\vartheta]$, see
Example \ref{hypersurface}). When $R$ is Golod with $\codepth
R\ge2$ Theorem \ref{structure2} does not apply, because $R$ is
not Gorenstein.  Neither does Theorem \ref{structure}(3):

 \begin{proposition}
 \label{golod}
If $R$ is a Golod ring with $\codepth R\ge2$, then
$\depth{}\CE=1$.
 \end{proposition}

 \begin{proof}
If $\cD$ is the universal enveloping algebra of $\pi^{\ges2}_R$, then
$\depth{}\cD\le1$, see \eqref{golod-def} and \eqref{tensor-algebra}, so
$\depth{}\CE\le1$ by Lemma \ref{subalgebra}; equality holds by
Lemma \ref{singular}.
  \end{proof}

Thus, new tools are needed to study stable cohomology over Golod
rings of higher codimension.  Local rings with radical square zero
are the simplest example of Golod rings.  They are the subject of the
next theorem, proved at the end of the section; its notation
and hypotheses are in force for the rest of the section.

\begin{theorem}
 \label{square-zero}
Let $(R,\fm,k)$ be a local ring with $\fm^2=0$ and $\edim R=e\ge2$.

The following assertions then hold.
 \begin{enumerate}[\rm\quad(1)]
 \item
The exact sequence of Theorem \emph{\ref{structure}(1)} has the
form
 \begin{equation*}
0\lra\CE\xra{\ \iota\ }\CS\xra{\ \eth\ } (\shift\CI)^{e}\oplus
\coprod_{i=0}^\infty(\shift^{-i}\CI)^{e^{i}(e^2-1)}\lra0
  \end{equation*}
 \item
The left $\CE$-submodule $\iota(\CE)$ has no direct
complement in $\CS$.
 \item
One has $\Ext n\CE k\CS=0$ for all $n\in\BZ$.
 \end{enumerate}
  \end{theorem}

For the next lemma we introduce shorthand
notation: When $\beta,\beta'\col M\to N$ are $R$-linear maps we write
$\beta'\equiv\beta\pmod\fm$ in place of $(\beta'-\beta)(M)\subseteq\fm N$.

\begin{lemma}
 \label{bar}
Let $\delta\col U\to R$ be the composition of a projective cover
$U\to\fm$ of the $R$-module $\fm$ with the inclusion $\fm\subseteq
R$. For $i\in\BZ$, let $\dd_{i+1}$ be the $R$-linear map
 \[
F_{i+1}=U^{\otimes(i+1)}=U\otimes_RU^{\otimes i}
\xra{\delta\otimes_R{U^{\otimes i}}}R\otimes_RU^{\otimes i}
=U^{\otimes i}=F_i
 \]
where $U^{\otimes i}$ denotes the $i$th tensor power of $U$ over
$R$, with the conventions $U^{\otimes 0}=R$ and $U^{\otimes i}=0$
for $i<0$. The following then hold.
\begin{enumerate}[\rm\quad(1)]
 \item
The pair $(F,\dd)$ is a minimal $R$-free resolution of $k$.
 \item
There are equalities $\dd(F_{i+1})=\fm F_i$ for all $i\in\BZ$.
 \item
When $h$ and $i$ are integers, such that $i\ge\max\{0,-h\}$, a
diagram
 \[
 \xymatrixrowsep{2.5pc}\xymatrixcolsep{1.3pc}
  \xymatrix{
U\otimes_RF_{i}
 \ar@{=}[r]
&F_{i+1}
 \ar@{->}[rr]^{\dd_{i+1}}
 \ar@{->}[d]_{\beta'}
&&F_{i}
 \ar@{->}[d]^{\beta}
  \\
U\otimes_RF_{i+h}
 \ar@{=}[r]
&F_{h+i+1}
 \ar@{->}[rr]^{\dd_{h+i+1}}
&&F_{h+i}
 }
 \]
of $R$-linear maps commutes if and only if
$\beta'\equiv U\otimes_R\beta\pmod\fm$.
 \end{enumerate}
 \end{lemma}

\begin{noremark}
When $R$ is a $k$-algebra, the minimal resolution $F$ described in Part
(1) of the theorem coincides with the bar construction of $R$ over $k$.
 \end{noremark}

\begin{proof}
(2) This follows directly from the definition of $\dd$.

(1) {}From $\fm^2=0$ one gets $\dd^2=0$. For $i\ge1$ set $B_{i-1}=
\dd_{i}(F_{i})$ and $Z_i=\Ker(\dd_i)$. {}From (2) one obtains
$\ell(B_{i-1}) =e^{i}$, where $\ell$ denotes length over $R$. The
exact sequence $0\to Z_{i}\to F_i\to B_{i-1}\to0$ yields
$\ell(Z_{i})=(e+1)e^i-e^i=e^{i+1}=\ell(B_i)$.  Thus, $\HH iF=0$
for $i\ge1$ and $\HH0F=k$, so $F$ is a minimal free resolution of
$k$.

(3)  Pick $u\otimes v\in U\otimes_RF_i=F_{i+1}$.  By definition, one has
\[
\beta\dd_{i+1}(u\otimes v)
=\beta(\delta(u)v)=\delta(u)\beta(v)=\dd_{h+i+1}(u\otimes\beta(v))\,.
\]
Therefore, an equality $\beta\dd_{i+1}(u\otimes v)=
\dd_{h+i+1}\beta'(u\otimes v)$ holds if and only if one has
\[
\beta'(u\otimes v)-u\otimes\beta(v)\in\Ker(\dd_{h+i+1})\,.
 \] 
As (1) and (2) yield $\Ker(\dd_{h+i+1})=\dd_{h+i+2}(F_{h+i+2})=\fm F_{h+i+1}$,
the inclusion above is equivalent to the relation
$\beta'\equiv U\otimes_R\beta\pmod\fm$.
  \end{proof}

\begin{lemma}
 \label{description}
Let $F$ denote the minimal free resolution from the preceding
lemma, let $\varkappa\in\Hom RFF_h$ be a homomorphism of
complexes, let $\wh\varkappa$ denote its image in $\Homv RFF_h$,
and set $m=\max\{0,-h\}$. The following then hold.
 \begin{enumerate}[\rm\quad(1)]
  \item
The map $\varkappa$ (respectively, $\wh\varkappa$) is a boundary if and
only if for $i=m$ (respectively, for some $i\ge m$) and for all
$j\ge0$ one has
 \[
\varkappa_{i+j}\equiv0\pmod\fm\,.
 \]
  \item
The map $\varkappa$ (respectively, $\wh\varkappa$) is a cycle if and only
if for $i=m$ (respectively, for some $i \ge m$) and for all
$j\ge0$ one has
 \[
\varkappa_{i+j}\equiv U^{\otimes j}\otimes_R(-1)^{hj}\varkappa_{i}\pmod\fm\,.
 \]
  \end{enumerate}
\end{lemma}

\begin{proof}
(1)  By definition, $\varkappa$ (respectively, $\wh\varkappa$) is a
boundary if and only if there exists a homomorphism $\chi\in\Hom RFF$
of degree $h+1$, such that an equality
 \begin{equation}
  \label{homotopic}
  \tag{\ref{description}.3}
\varkappa_{i+j}=\dd_{h+i+j+1}\chi_{i+j}+(-1)^{h}\chi_{i+j-1}\dd_{i+j}
 \end{equation}
holds for $i=m$ (respectively, for some $i\ge m$) and for all
$j\ge 0$.

If $\chi$ exists, then $\varkappa_{i+j}\equiv0\pmod\fm$ for $i,j$ as
above, because $\dd\equiv0\pmod\fm$.

Conversely, assume $\varkappa_{i+j}\equiv0\pmod\fm$ holds for $i,j$
as above.  We construct $\chi_{i+j}$ by induction on $j$.  Setting
$\chi_{i+j}=0$ for $j<0$, we may assume that $\chi_{i+j-1}$ has been
defined for some $j\ge 0$.  One then has the relations
 \begin{align*}
\big(\varkappa_{i+j}-(-1)^{h}\chi_{i+j-1}\dd_{i+j}\big)(F_{i+j})
&\subseteq\fm F_{h+i+j}+\chi_{i+j-1}(\fm F_{i+j-1})\\
&\subseteq \fm F_{h+i+j}\\
&=\dd_{h+i+j+1}(F_{h+i+j+1})\,.
 \end{align*}
with equality given by Lemma \ref{bar}(2). Since $F_{i+j}$ is free,
one can find a homomorphism $\chi_{i+j}\col F_{i+j}\to F_{h+i+j+1}$
satisfying equation \eqref{homotopic}.

(2) By definition, $\varkappa$ (respectively, $\wh\varkappa$)
is a cycle if and only if $\dd_{h+i+j}\varkappa_{i+j}=
(-1)^h\varkappa_{i+j-1}\dd_{i+j}$ holds for $i=m$ (respectively, for some
$i\ge m$) and for all $j\ge 0$. Iterated applications of (\ref{bar}.3)
yield the desired assertion.
  \end{proof}

It is convenient to describe $\CE$ and $\CS$ as subrings of
infinite matrix rings.

\begin{chunk}
 \label{matrices}
For every pair $(m,r)\in\BZ\times\BN$, satisfying $m+r\ge0$, let
$\mat{e^r\times e^{r+m}}k$ denote the $k$-vector space of $e^r\times
e^{r+m}$ matrices with elements in $k$.  Let $\mat{\infty}k$ be the
$k$-algebra of all row-and-column-finite matrices with elements in
$k$, under ordinary matrix product. For every matrix $C=(c_{ij})\in
\mat{e^r\times e^{r+m}}k$ let $C_\infty=(c^\infty_{pq})\in
\mat{\infty}k$ be the matrix with blocks $C$ along a line of slope $e^{-m}$:
 \[
c^\infty_{pq}=\begin{cases}
 c_{ij}&\text{if }p=i+le^{r}\text{ and }
  q=j+le^{r+m}\text{ for some }l\ge0\,;
 \\
0&\text{ otherwise}\,.
 \end{cases} \]

It is clear that for each $m\in\BZ$ the following subset of
$\mat{\infty}k$ is a $k$-subspace:
 \[
\CC^m=
 \left\{\,C_\infty\in \mat{\infty}k
\,\left|\,
\begin{gathered}{C}\in\mat{e^r\times e^{r+m}}k\text{ for some}
 \\
\text{pair $(m,r)$ with }m+r\ge0
 \end{gathered}\,
\right\}\right.\,.
 \]
A key observation is that, furthermore, the following relations hold:
 \begin{alignat*}{2}
\CC^m\cdot\CC^n&\subseteq\CC^{m+n}&\quad&\text{for all}\quad
m,n\in\BZ\,;
 \\
\CC^m\cap\CC^n&=0&\quad&\text{when}\quad m\ne n\,.
 \end{alignat*}

Indeed, let $C_\infty\in\CC^m$ be as above, and let $D_\infty\in\CC^n$
be obtained from a matrix $D\in\mat{e^s\times e^{s+n}}k$,
where $(n,s)\in\BZ\times\BN$ satisfy $n+s\geq 0$.  One then has
$C_\infty=C'_\infty$, where $C'\in\mat{e^{r+s}\times e^{r+s+m}}k$ is
the block diagonal matrix with $e^s$ copies of $C$ along a line of slope
$e^{-m}$; also, $D_\infty=D'_\infty$ where $D'\in\mat{e^{r+s+m}\times
e^{r+s+m+n}}k$ is the block diagonal matrix with $e^{r+m}$ copies of $D$
along a line of slope $e^{-n}$.  Thus, we get
  \[
C'D'\in\mat{e^{r+s}\times e^{r+s+m+n}}k
  \quad\text{and}\quad
C_\infty\cdot D_\infty=(C'D')_\infty\in\CC^{m+n}\,.
  \]
This proves the inclusion.  For the equality, assume $C_\infty=D_\infty\in
\CC^m\cap\CC^n$ with $m\not=n$.  Since the lines with slopes $e^{-m}$
and $e^{-n}$ diverge, for $l\gg0$ the blocks $C$ forming the matrix
$C_\infty$ are entirely contained in an area of the matrix $D_\infty$
where every element is equal to $0$.  Thus, one has $C=0$, and
consequently $C_\infty=0$.

By the discussion above, $\CC=(\CC^m)_{m\in\BZ}$ is a graded
$k$-algebra with unit $1_\infty$.  

We define a graded subalgebra $\CA$ of $\CC$ as follows. Set
 \begin{align*}
A^{(m)}_{i}&=(E^{(m)}_i)_\infty\in\CC^m
\quad\text{for each pair $(m,i)\in\BN\times[1,e]$, where}
 \\
E^{(m)}_i&=[0,\dots,0,1,0,\dots,0]\in\mat{1\times e^m}k
\quad\text{with $1$ in $i$th position}\,.
 \end{align*}
Let $\CA$ be the subalgebra generated over
$k1_\infty$ by $A^{(1)}_{1},\dots,A^{(1)}_{e}\in\CC^1$.  It is easy
to see:
 \[
A^{(m)}_{i}\,A^{(n)}_{j}=A^{(m+n)}_{(i-1)e^n+j}\,.
 \]
Thus, $\CA$ has a $k$-basis consisting of matrices with exactly one
non-zero entry in every row, and in distinct basis elements this entry
occurs in distinct columns.
 \end{chunk}

Part (2) of the next lemma is known:  For a $k$-algebra $R$ it is obtained
by computing $\Ext{}Rkk$ as the cohomology of the cobar construction,
which in this case is the tensor algebra on $\Hom k{\fm/\fm^2}k$ with
zero differential; the general case can be found in \cite[Thm.\ 1,
Cor.\ 3]{Ro}.  A proof is included for completeness.

\begin{lemma}
 \label{matrix-algebra}
In the notation of \emph{\eqref{matrices}} the following hold.
 \begin{enumerate}[\rm\quad(1)]
 \item
There is an isomorphism of graded $k$-algebras $\CS\cong\CC$,
inducing $\iota(\CE)\cong\CA$.
 \item The associative $k$-algebra $\CE$ is freely generated
over $k$ by a $k$-basis of $\CE^1$.
 \item
One has $\Ext n\CE k\CE=0$ if $n\ne1$ and
 \[
\rank_k\Ext1\CE k\CE^i
 =\begin{cases}
0\quad&\text{for}\quad i\le-2\,;
 \\
e&\text{for}\quad i=-1\,;
 \\
e^{i}(e^2-1)&\text{for}\quad i\ge0\,.
 \end{cases}
 \]
 \end{enumerate}
  \end{lemma}

\begin{proof}
Let $F\to k$ be the minimal free resolution from Lemma
\ref{bar}. We fix a basis $X_1=\{x_1,\dots,x_e\}$ of $U$ over
$R$. For each $i\ge0$ it canonically provides a basis $X_i$ of
$F_i=U^{\otimes i}$ over $R$, and thus a basis $\ov X_i$ of
$F_i/\fm F_i$ over $k$.

(1)  Each $\sigma\in\CS^h$ is the class of a cycle
$\wh\varkappa\in\Homv RFF_{-h}$, where $\varkappa\col F\to F$ is a
homomorphism of complexes of $R$-modules of degree $-h$.  Thus,
$\varkappa_{\ges s}\col F_{\ges s}\to F_{\ges s-h}$ is a chain map
for some integer $s$. It induces a $k$-linear map
 \[
\HH s\varkappa=\ov\varkappa_s\col U^{\otimes(s+h)}/\fm
U^{\otimes(s+h)}\lra U^{\otimes s}/\fm U^{\otimes s}\,.
 \]
Let $S\in\mat{e^{s}\times e^{s+h}}k$ be the matrix of
$\ov\varkappa_s$ in the bases $\ov X_{s+h}$ and $\ov X_{s}$, and
form the matrix $S_\infty\in\CC^h$. Lemma \ref{description}
shows that $S_\infty$ does not depend on the choices of
$\varkappa$ or $s$, so setting $\alpha(\sigma)=S_\infty$ one
obtains a map $\alpha\col\CS\to\CC$. The definitions of the
products in $\CS$ and $\CC$ imply that $\alpha$ is a homomorphism
of algebras. Part (1) of Lemma \ref{description} shows that
$\alpha$ is injective and Part (2) that it is surjective.

Let $\{\xi_1,\dots,\xi_e\}$ be the basis of $\CE^1$ dual to the
basis $\ov X_1$ of $U/\fm U$.  By definition,
$\alpha(\iota(\xi_j))=A^{(1)}_{j}$ for $j=1,\dots,e$, so $\alpha$
maps the subalgebra of $\CE$ generated by $\{\xi_1,\dots,\xi_e\}$
surjectively onto the subalgebra $\CA$ of $\CC$. Thus, for each
$i\in\BZ$ one has
 \[
e^i=\rank_k\CE^i\ge\rank_k\CA^i=e^i\,,
 \]
with equalities given by the constructions in Lemma \ref{bar}
and in \eqref{matrices}, respectively.  Thus, $\alpha$ restricts
to an isomorphism $\iota(\CE)\cong\CA$.

(2) As shown above, $\rank_k\CE^i=e^i$ and $\CE$ is generated by
$e$ elements of degree $1$.  It follows that there are no
relations between the generators of $\CE$.

(3) In view of (2), the exact sequence of \eqref{tensor-algebra}
yields $\Ext n\CE k\CE=0$ for $n\ne0,1$ and an exact sequence of
graded vector spaces
 \[
0\lra\Hom\CE k\CE\lra\CE\xra{\ \dd^*\ }\shift\CE^e\lra\Ext1\CE
k\CE\lra0
 \]
where $\dd^*(\rho)=(\xi_1\rho,\dots,\xi_e\rho)$.  We see that
$\Hom\CE k\CE=0$.  Counting $k$-ranks in the sequence above we now
get the desired expressions for $\Ext1\CE k\CE^i$.
   \end{proof}

 \label{square-zero:proof}
\begin{proof}[Proof of Theorem \emph{\ref{square-zero}}]
The notation from the preceding proof stays in force.

(1) Theorem \ref{structure}(1) yields an exact sequence of graded
left $\CE$-modules
 \begin{equation}
\label{old} 0\lra\CE\xra{\ \iota\ }\CS\xra{\ \eth\ }
\coprod_{i=-1}^\infty(\shift^{-i}\CI)^{\mu^{i+1}}\lra0
 \end{equation}
where $\CI=\Hom k\CE k$ and $\mu^i=\rank_k\Ext iRkR$.  These
numbers are given by:
\begin{equation}
 \label{moo}
\mu^i
 =\begin{cases}
0\quad&\text{for}\quad i\le-1\,;
 \\
e&\text{for}\quad i=0\,;
 \\
e^{i-1}(e^2-1)&\text{for}\quad i\ge1\,.
 \end{cases}
 \end{equation}

Indeed, the hypothesis $\fm^2=0$ yields an exact sequence of
$R$-modules
 \[
0\lra k^e\lra R\lra k\lra0
 \]
It induces an exact sequence of homomorphisms of $k$-vector spaces
 \[
0\lra\Hom RkR\lra R\lra\Hom R{k^e}R\lra\Ext1RkR\lra0
\]
giving the values of $\mu^i$ in \eqref{moo} for $i\le1$.  The same
sequence yields isomorphisms
 \[
\Ext{i-1}R{k^e}R\cong\Ext{i}RkR \quad\text{for all}\quad i\ge2\,,
 \]
which imply $\mu^{i}=e\mu^{i-1}$ for all $i\ge2$; the last equality
in \eqref{moo} follows.

(3) The $k$-algebra $\CE$ is free by Lemma \ref{matrix-algebra}(2),
so the resolution of $k$ displayed in \eqref{tensor-algebra} gives 
$\Ext n\CE k\CS=0$ for all $n\ne0,1$.

Next we prove $\Hom\CE k\CS=0$.  By \eqref{depth0}, this is
equivalent to the following assertion:  If $\varkappa\in\Hom
RFF_h$ is such that $\wh\varkappa\in\Homv RFF_h$ is a cycle and
$\cls{\wh\varkappa}\in\Extv{-h}Rkk$ satisfies $\Ext{\ges
s}Rkk\cdot\cls{\wh\varkappa}=0$ for some $s\ge0$, then
$\cls{\wh\varkappa}=0$.

Set $m=\max\{0,-h\}$.  As $\wh\varkappa$ is a cycle, for some integer
$i$ with $i\ge m$ one has
 \[
\dd_{h+i+j}\varkappa_{i+j}=(-1)^{h}\varkappa_{i+j-1}\dd_{i+j}
\quad\text{for all}\quad j\ge0\,.
 \]
After increasing $i$ or $s$ (if necessary) we may assume
$s=i+h\ge0$.

Let $\xi\in\Hom RFF^s$ be a chain map.  As
$\cls\xi\cdot\cls{\wh\varkappa}=0$, Lemma \ref{description}(1)
yields
 \[
(\xi\varkappa)_{i+j}(F_{i+j})\subseteq\fm F_{j} \quad\text{for
all}\quad j\gg0\,.
 \]
On the other hand, one has $ \dd_{j}\circ(\xi\varkappa)_{i+j}
=(-1)^{h+s}(\xi\varkappa)_{i+j-1}\circ\dd_{i+j}$ for all $j\ge0$,
so Lemma \ref{bar}(3) implies that these inclusions hold, in
fact, for all $j\ge0$.

Recall that $X_i$ denotes the standard basis of $F_i=U^{\otimes
i}$. For each $u\in X_{i}$ one has $\varkappa_{i}(u)=\sum_{v\in
X_{i+h}}a_{uv}v$ with uniquely defined $a_{uv}\in R$.  On the
other hand, for each $v\in X_{h+i}$ the $R$-linear map $F_{h+i}\to
F_0=R$ sending $v$ to $1$ and every $v'\in
X_{h+i}\smallsetminus\{v\}$ to $0$ extends to a chain map
$\xi_v\col F\to F$ of degree $-h-i$. We get
 \[
a_{uv}=(\xi_v)_{h+i}\bigg(\sum_{v\in X_{h+i}}a_{uv}v\bigg)=
(\xi_v)_{h+i}\big(\varkappa_{i}(u)\big)=(\xi_v\varkappa)_{i}(u)\in\fm
 \]
hence $\varkappa_{i}(F_{i})\subseteq\fm F_{h+i}$ holds.  Lemma
\ref{bar}(3) now yields $\varkappa_{i}(F_{i}) \subseteq\fm
F_{h+i}$ for all $i\ge0$, from where we conclude
$\cls{\wh\varkappa}=0$, see Lemma \ref{description}(1).

Finally, we prove $\Ext1\CE k\CS=0$.  As $\Hom\CE k\CS=0$, the
exact sequence \eqref{old} of graded left $\CE$-modules induces an
exact sequence of graded vector spaces
 \begin{equation}
\label{useful}
\begin{gathered}
\xymatrixrowsep{.1pc} \xymatrixcolsep{1.5pc}\xymatrix{
0\ar@{->}[r]
&\BigHom\CE{k}{{\ds\coprod_{i=-1}^\infty(\shift^{-i}\CI)^{\mu^{i+1}}}}
\ar@{->}[r]
&\Ext1\CE k\CE \ar@{->}[r]^-{\iota_*}
&\Ext1\CE k\CS
\\
{\hphantom{0}}\ar@{->}[r]^-{\eth_*}
&\BigExt1\CE{k}{\ds\coprod_{i=-1}^\infty(\shift^{-i}\CI)^{\mu^{i+1}}}
}
 \end{gathered}
 \end{equation}
Since $k$ has a resolution by free $\CE$-modules of finite rank,
see \eqref{tensor-algebra}, the functors $\Ext n{\CE}k-$ commute with
direct sums.  The graded left $\CE$-module $\CI=\Hom k\CE k$
satisfies $\Hom\CE k\CI\cong k$ and is injective, so we obtain
 \begin{align}
  \label{coincidence}
\BigHom\CE{k}{\coprod_{i=-1}^\infty(\shift^{-i}\CI)^{\mu^{i+1}}}
&\cong\coprod_{i=-1}^\infty\shift^{-i}k^{\mu^{i+1}}\,;
 \\
 \label{nothing}
\BigExt1\CE
k{\coprod_{i=-1}^\infty(\shift^{-i}\CI)^{\mu^{i+1}}}&=0\,.
 \end{align}
Comparing \eqref{coincidence}, \eqref{moo}, and Lemma
\ref{matrix-algebra}(3) we see that in \eqref{useful} the map
$\iota_*$ is bijective. In view of \eqref{nothing} this implies
$\Ext1\CE k\CS=0$.

(2) If $\CS=\iota(\CE)\oplus\CT'$ for some left graded $\CE$-submodule
$\CT'$ of $\CS$, then Theorem \ref{structure}(2) yields $\CT'=\gam{\CS}$,
hence $\gam{\CS}^i=\CS^i$ for $i<0$.  One has $\CS^i\ne0$, see
Corollary \ref{regular}, so \eqref{linked} and \eqref{depth0} imply
$\depth\CE\CS=\depth{\CE}{\gam{\CS}}=0$.  This contradict (3).
   \end{proof}

We finish by applying lemmas used in the proof of Theorem
\ref{square-zero} to show that the action of absolute cohomology on
bounded cohomology from the \emph{right} may be far from nilpotent---in
contrast to the action from the \emph{left}, cf.\ Lemma \ref{nilpotent}.

 \begin{example}
 \label{non-nilpotent}
If $\fm^2=0$ and $\edim R=e\ge2$, then for every $n<0$ there
exists $\beta_n\in\Extb{n}Rkk$, such that 
$\{\varepsilon\in\Ext{}Rkk\var\beta_n\cdot\varepsilon=0\}$
is equal to $0$.

By Theorem \ref{structure}(1) and Lemma \ref{matrix-algebra}(1), 
it suffices to show that each matrix 
 \[
C_n=(C^{(n)})_\infty\in\CC^n\,,
 \quad\text{where $n<0$ and}\quad
C^{(n)}=\begin{bmatrix}1\\0\\
\vdots\\0\end{bmatrix}\in\mat{e^{-n}\times 1}k\,,
  \]
satisfies $(C_n\CA)\cap\CA=0$.   Indeed, for each $A\in\mat{\infty}k$ the 
matrix $C_nA$ is obtained by inserting $(e^{-n}-1)$ 
rows of zeroes between every pair of adjacent rows of $A$.  On the
other hand, the $k$-basis of $\CA$ described in \eqref{matrices}
shows that each non-zero matrix in $\CA$ has a non-zero entry in 
every row.
 \end{example}

\appendix

\section{Depth over graded algebras}
\label{Depth over graded algebras}

In this appendix $k$ is a field, $\CA$ is a graded $k$-algebra with
$\CA^0=k$ and $\CA^i=0$ for all $i<0$.  Throughout, $\CM$ denotes a
graded left $\CA$-module.  By a customary abuse of notation, we let $k$
denote also the graded $\CA$-module $\CA/\CA^{\ges1}$.

\begin{chunk}
 \label{depth}
The \emph{depth} of $\CM$ over $\CA$ is the number
\[
\depth{\CA}\CM=\inf\{n\in\BN\mid\Ext n\CA k\CM\ne0\}\,,
\]
see \cite{FHJLT}.  Clearly, one has $0\le\depth\CA\CM\le\infty$,
and $\depth\CA\CM=\infty$ holds if and only if $\Ext {}\CA k\CM=0$.
We systematically write $\depth\CA\CA$ in place of $\depth{}\CA$.
 \end{chunk}

Here we collect general facts about depth, for use in the body of the paper.

The long exact sequence of functors $\Ext n{\CA}k-$ yields a familiar
formula:

\begin{chunk}
 \label{linked}
The depths of the graded left $\CA$-modules appearing in an exact
sequence $0\to\CL\to\CM\to\CN\to0$ are linked by an inequality
 \[
\depth\CA\CM\ge\inf\big\{\depth\CA\CL\,,\depth\CA\CN\,\big\}\,.
 \]
Equality holds when the sequence splits, or when $\depth\CA\CL
\ne\depth\CA\CN+1$.
 \end{chunk}

For finite modules over finitely generated commutative algebras
depth measures lengths of regular sequences.  In general, only a
weaker statement holds.

\begin{chunk}
 \label{rees}
Assume $\CA=\CA'/(\vartheta')$ and $\CM=\CM'/\vartheta'\CM'$ for
some graded $k$-algebra $\CA'$, a central element
$\vartheta'\in\CA'^{\ges1}$, and a graded left $\CA'$-module
$\CM'$.

When $\vartheta'$ is a non-zero-divisor on $\CM'$ the following
hold.
 \begin{enumerate}[\rm\quad(1)]
  \item
$\depth{\CA'}{\CM'}=\depth{\CA'}{\CM}+1$.
 \item
If $\vartheta'$ is a non-zero-divisor on $\CA'$, then
$\depth{\CA}{\CM}=\depth{\CA'}{\CM'}-1$.
 \end{enumerate}

Indeed, for each $n\ge0$ there are isomorphisms of graded
$k$-vector spaces
\begin{gather*}
\Ext{n-1}{\CA'}k{\CM}\cong\Ext{n-1}{\CA'}k{\CM'}\oplus\Ext{n}{\CA'}k{\CM'}\\
\Ext{n-1}{\CA}k\CM\cong\Ext{n}{\CA'}k{\CM'}
\end{gather*}
obtained by transcribing Rees' classical argument for commutative
algebras.
 \end{chunk}

Pursuing the analogy with commutative algebra, we define
(left) section functors.

\begin{chunk}
 \label{gamma}
For each $i\ge0$ the graded subspace $\CA^{\ges i}$ of $\CA$ is a
two-sided ideal, so the following subspaces of $\CM$ are graded
left $\CA$-submodules:
 \setcounter{equation}{2}
\begin{equation*}
\gami i\CM=\big\{\mu\in\CM\ \big|\ \CA^{\ges i}\cdot\mu=0\}
 \quad\text{and}\quad
\gam\CM=\bigcup_{i=0}^\infty\gami i\CM\,.
 \end{equation*}
   \end{chunk}

Section functors carry information on the vanishing of depth.

\begin{chunk}
\label{depth0}
$\depth\CA\CM=0$ if and only if $\gami 1\CM\ne0$, if and only if
$\gam\CM\ne0$.

Indeed, the equivalence of the first two conditions comes from the
isomorphism
\[
\Hom\CA k\CM\cong\{\mu\in\CM\mid\CA^{\ges1}\cdot\mu=0\}=\gami1M\,.
\]
It is clear that the second condition implies the third one.
Conversely, if $\gam\CM\ne0$, then $\CA^{\ges i}\cdot\mu=0$ for
some $\mu\in\CM\smallsetminus\{0\}$ and some integer $i\ge1$.  Choosing
$i$ minimal with this property, for $\CN=\CA^{\ges i-1}\cdot\mu$ we
get $0\ne\CN\subseteq\gami1\CM\ne0$.
   \end{chunk}

Several applications of depth in the body of the paper hinge on
the next result.

\begin{proposition}
 \label{torsion}
If a graded left $\CA$-module $\CK$ satisfies $\CK=\gam\CK$, then
 \[
\depth\CA\CM\le\inf\{n\in\BN\mid\Ext n\CA\CK\CM\ne0\}\,.
 \]
Equality holds if $\CK=\gami i\CK\ne0$ for some integer $i\ge1$.
 \end{proposition}

 \begin{proof}
Set $m=\depth{\CA}{\CM}$.  We prove the last assertion by induction
on $i$.  If $i=1$, then $\CK$ is a direct sum of shifts of copies of $k$,
so $\Ext n{\CA}\CK\CM$ is a direct product of shifts of $\Ext n{\CA}k\CM$,
and the assertion is clear.  When $i>1$ the exact sequence
 \[
0\lra\gami{i-1}\CK\lra\gami{i}\CK\lra\CL\lra0
 \]
of graded left $\CA$-modules, where $\CL=\gami{i}\CK/\gami{i-1}\CK$,
induces an exact sequence
 \[
\Ext{n-1}{\CA}{\gami{i-1}\CK}{\CM}
 \lra\Ext n{\CA}{\CL}{\CM}
\lra\Ext n{\CA}{\gami{i}\CK}{\CM}
 \lra\Ext n{\CA}{\gami{i-1}\CK}{\CM} \]
for each $n$.  By the base of the induction, $\Ext n{\CA}{\CL}{\CM}$
vanishes for $n<m$ and does not for $n=m$.  For $n<m$ the first
and last terms vanish by the induction assumption. Thus, $\Ext
n{\CA}{\gami{i}\CK}{\CM}$ vanishes for $n<m$ and does not for $n=m$.

In general, each $\gami i\CK$ is a graded left submodule of $\CK$,
so there are exact sequences
\[
0\lra\varprojlim_i{\!}^{1}\Ext{n-1}{\CA}{\gami i\CK}{\CM}
 \lra\Ext n{\CA}{\CK}{\CM}
\lra\varprojlim_i\Ext n{\CA}{\gami i\CK}{\CM}\lra 0
\]
for all $n\ge0$, where $\varprojlim{\!}^{1}$ is the first right derived
functor of $\varprojlim$, see \cite[(3.5.10)]{We}. In view
of the finite case, these sequences yield $\Ext n{\CA}{\CK}{\CM}=0$
for $n<m$.
 \end{proof}

\begin{corollary}
 \label{extension1}
If $\CB$ is a graded $k$-subalgebra of $\CA$, such that $\CA$ is
free as a graded right $\CB$-module and one has $\CA^{\ges
i}\subseteq \CA\cdot\CB^{\ges1}$ for some $i\ges1$, then
 \[
\depth{\CB}\CM=\depth{\CA}\CM\,.
 \]
 \end{corollary}

\begin{proof}
For $\ov\CA=\CA/(\CA\cdot\CB^{\ges1})$ standard arguments yield
$\Ext{n}\CA{\ov\CA}\CM\cong\Ext n\CB k\CM$ for each $n\in\BZ$.
Since $\ov\CA=\gami i{\ov\CA}$, the proposition yields the desired
equality.
 \end{proof}

A graded $k$-subalgebra $\CB$ of $\CA$ is said to be \emph{normal}
if $\CB^{\ges1}\cdot\CA=\CA\cdot\CB^{\ges1}$.

\begin{corollary}
 \label{extension2}
Assume $\rank_k\CA^i$ is finite for each $i$ and $\CB$ is a
normal subalgebra of $\CA$.  A finite subset $E\subseteq\CA$ is a
basis of $\CA$ as a graded right $\CB$-module if and only if it is
a basis of $\CA$ as a graded left $\CB$-module.  When such a set $E$
exists, one has
 \[
\depth{}\CB=\depth{}\CA\,.
\]
 \end{corollary}

\begin{proof}
By symmetry, to prove the first assertion it suffices to show that
if $E$ is a basis of $\CA$ as a graded right $\CB$-module, then it
is one as a left $\CB$-module. The image $\ov E$ of $E$ in
$\ov\CA=\CA/\CA\cdot\CB^{\ges1}$ is a $k$-basis of $\ov\CA$.  The
map
 \[
\gamma\col\ov\CA\otimes_k\CB\lra\CA\quad\text{given by}\quad
\gamma\bigg(\sum_{e\in E}a_e\ov e\otimes\delta_e\bigg)=\sum_{e\in
E}(-1)^{|e||\delta_e|}a_e\delta_e e
 \]
is a morphism of graded left $\CB$-modules. Since $\CB$ is normal
in $\CA$ one has
 \[
k\otimes_\CB\CA\cong\CA/(\CB^{\ges1}\cdot\CA)=\CA/(\CA\cdot\CB^{\ges1})
=\ov\CA\,,
 \]
so $k\otimes_\CB\gamma$ is bijective.  By (a graded version of)
Nakayama's Lemma the map $\gamma$ is then surjective. Comparison
of $k$-ranks shows that it is bijective.  Thus, $E$ is a basis of
$\CA$ as left $\CB$-module.  When a basis $E$ as above exists one
has $\rank_k\ov\CA<\infty$, whence the first isomorphism below;
the isomorphism $\gamma$ induces the second one:
 \[
\ov\CA\otimes_k\Ext{n}\CB k\CB\cong\Ext n\CB
k{\ov\CA\otimes_k\CB}\cong\Ext n\CB k\CA\,.
 \]
They yield $\depth{}\CB=\depth{\CB}\CA$, and Corollary
\ref{extension1} gives $\depth{\CB}\CA =\depth{}\CA$.
 \end{proof}


\section*{Acknowledgments}

The authors thank Dave Benson, Ragnar Buchweitz, Lars Christensen,
Yves F\'elix, Srikanth Iyengar, Frank Moore, R\u{a}zvan Veliche, Dan
Zacharia, and the referee for useful comments at various stages of the
long evolution of this paper.


\end{document}